\documentclass[a4paper,11pt,fleqn]{article}
\usepackage{amsmath,amssymb,amsthm,graphicx,subfigure,float,caption,epstopdf,tabularx,color, bm,amsfonts,epic}
\usepackage[top=1in, bottom=1in, left=1.25in, right=1.25in]{geometry}
\usepackage{appendix}
\allowdisplaybreaks
\newtheorem{thm}{Theorem}[section]
\newtheorem{lem}{Lemma}[section]
\newtheorem{coro}{Corollary}[section]
\newtheorem{den}{Definition}[section]

\newtheorem{rmk}{Remark}[section]
\newtheorem*{prf}{Proof}
\numberwithin{equation}{section}

\usepackage{indentfirst}
\graphicspath{{Fig}}

\begin{document}
\title{A novel linearized and momentum-preserving Fourier pseudo-spectral scheme for the Rosenau-Korteweg de Vries equation}
\author{Chaolong Jiang$^1$,\ Jin Cui$^{2,3}$,\ Wenjun Cai$^2$,\ Yushun Wang$^2$\footnote{Correspondence author. Email:
wangyushun@njnu.edu.cn.}\\
{\small $^1$ School of Statistics and Mathematics, }\\
{\small Yunnan University of Finance and Economics, Kunming 650221, China}\\
{\small $^2$ Jiangsu Provincial Key Laboratory for NSLSCS,}\\
{\small School of Mathematical Sciences,  Nanjing Normal University,}\\
{\small  Nanjing 210023, China}\\
{\small $^3$ Department of Basic Sciences,}\\
{\small Nanjing Vocational College of Information Technology,}\\
{\small Nanjing 210023, China}
}
\date{}
\maketitle

\begin{abstract}
In this paper, we design a novel linearized and momentum-preserving Fourier pseudo-spectral scheme to solve
the Rosenau-Korteweg de Vries equation.
With the aid of a new semi-norm equivalence between the Fourier pseudo-spectral method and the finite difference method, a prior bound of the numerical solution in discrete $L^{\infty}$-norm
is obtained from the discrete momentum conservation law. Subsequently, based on the energy method and the bound of the numerical solution, we show that, without any restriction on the mesh ratio, the scheme is convergent with order $O(N^{-s}+\tau^2)$  in discrete $L^\infty$-norm, where $N$ is
the number of collocation points used in the spectral method and $\tau$ is the time step. Numerical results are addressed to confirm our theoretical analysis. \\[2ex]
\end{abstract}

\section{Introduction}
To describe the dynamics of dense discrete systems, Rosenau \cite{Rosenau86,Rosenau88} derived the so-called Rosenau equation, as follows:
\begin{align}\label{R-equation}
u_t+u_{xxxxt}+u_x+uu_x=0,
\end{align}
where the existence and the uniqueness of the solution for \eqref{R-equation} were
proved by Park \cite{Park92}. On the other hand, for the further consideration
of the nonlinear wave, the viscous term $+u_{xxx}$ needs to be
included \cite{zuo09}
\begin{align}\label{R-equation:1.3}
u_t+u_{xxxxt}+u_{xxx}+u_x+uu_x=0.
\end{align}
%This equation is usually called the Rosenau-KdV equation.
%Rosenau \cite{Rosenau86,Rosenau88} derived the Rosenau equation to describe the dynamics of dense discrete systems.
%The existence and uniqueness of the solution for the Rosenau equation were discussed by Park \cite{Park92}.
%On the other hand, for the further consideration of the nonlinear wave, Zuo \cite{zuo09} added the viscous term $+u_{xxx}$ to the Rosenau equation.
Equation \eqref{R-equation:1.3} is usually called the Rosenau-KdV equation and the generalized case reads \cite{Esfahani11}
\begin{align}\label{GR-KDV1d}
 u_t+u_{xxxxt}+u_{xxx}+u_x+(u^p)_x=0,
\end{align}
where $1\le p<+\infty$ is a given integer.

The Rosenau-KdV equation has been studied theoretically and numerically in the literature. For the existence and uniqueness of the solution of the Rosenau-KdV equation, please refer to Refs. \cite{Esfahani11,zuo09}. Known strategies to solve the Rosenau-KdV equation numerically include finite difference methods \cite{AO15,CLZ18,HXH13,OAAK15}, a Fourier pseudo-spectral method \cite{CSW15}, etc. However, there are few works on the Rosenau-KdV equation in high dimensions. Thus, in this paper, we focus on developing a numerical method for the following generalized Rosenau-KdV (GR-KdV) equation in two dimensions (2D) \cite{AO15}
\begin{align}\label{GR-KDV_eq:1.1}
 u_t+\Delta^2u_t+\Delta u_x+(1+u^p)\mathbb{L}u=0,\ \ (x,y)\in\Omega,\ 0<t\leq T,
\end{align}
with the $(l_{1},l_{2})$-periodic boundary conditions
\begin{align*}%\label{eq:1.2}
&u(x,y,t)=u(x+l_{1},y,t),\ u(x,y,t)=u(x,y+l_{2},t),\ (x,y)\in\Omega,\ 0<t\leq T,
\end{align*}
and the initial condition
\begin{align*}%\label{eq:1.3}
u(x,y,0)=u_{0}(x,y),\ (x,y)\in\Omega,
\end{align*}
where $\Delta$ is the usual Laplace perator, $\mathbb{L}=\partial_x+\partial_y$, $\Omega=[x_L,x_R]\times[y_L,y_{R}]\subset\mathbb{R}^{2}$, $l_1=x_R-x_L$, $l_2=y_R-y_L$,
 and $u_0:=u_{0}(x,y)$ is a given real-valued initial data.
%When $\gamma=0$, the DNLS equation \eqref{3DNLD_eq:1.1} reduces to the classical NLS equation.
%Miles \cite{Miles67} reviewed and derived analytic formulae for $\gamma$, based on various kinds of physical dissipation that
%affect waves in deep water.
Under the periodic boundary conditions, the system \eqref{GR-KDV_eq:1.1} has the following %mass conservation law
%\begin{align*}%\label{eq:1.6}
%M(t)=\int_{\Omega}u d{\bm x}\equiv M(0),
%\end{align*}
% and
momentum conservation law
 \begin{align}\label{MCL}
 \mathcal{P}(t)=\iint_{\Omega}(u^2+|\Delta u|^2)dxdy\equiv \mathcal{P}(0),\ 0\le t\le T.
 \end{align}

%It is well-known that energy conservative law is an essential property in the theory of solitons, especially
%in an elastic collision of solitons.
 %Based on the basic rule that numerical methods should preserve the intrinsic
%properties of the original problems as much as possible, Feng \cite{feng86} first presented the
%concept of symplectic schemes for Hamiltonian systems and furthered the structure-preserving
%methods for the general conservative dynamical systems.
In recent years, there has been growing interest in geometric methods or structure-preserving methods, which can preserve as much as possible the intrinsic properties of the given dynamical system. It has been shown that, compared with traditional numerical methods, structure-preserving methods have excellent stability and superior performance in long time simulations. For more details, please refer to Refs. \cite{BR01,FQ10,ELW06} and references therein.
With the aid of the variational formulation \cite{marsden}, Cai et al. first derived some multi-symplectic schemes for the Rosenau-type equation \cite{CSW15}. More recently,
based on the multi-symplectic Hamiltonian formula \cite{BR01}, a new multi-symplectic scheme has been proposed for the Rosenau-type equation with the power law nonlinearity in Ref. \cite{CLZ18}. Besides the multi-symplectic structure, the
Rosenau-KdV equation also admits some invariants, such as the momentum conservation law \eqref{MCL}. In many significant cases, the ability to preserve
some invariant properties of the original differential equation is a criterion to judge the success
of a numerical simulation.
Thus, when discretizing such a conservative system in space and time, it is a natural idea to design numerical schemes that preserve rigorously
a discrete invariant. In Ref. \cite{HXH13}, the authors proposed a three-level linear finite difference scheme, which can preserve the conservation law \eqref{MCL}, for the Rosenau-KdV equation in 1D. In Ref. \cite{AO15}, Atouani and Omrani constructed two conservative schemes for the 2D case. However, most of the existing conservative schemes have only second order accuracy in space. To construct high order schemes in space, the Fourier spectral method \cite{ST06} was employed to discrete the Rosenau-KdV and two high order structure-preserving schemes were constructed in Ref. \cite{CLZ18}. However, the resulting schemes are fully implicit, which implies that one has to solve a system
  of nonlinear equations, at each time step. Thus, the first purpose of this paper is to develop a novel linearized and structure-preserving Fourier pseudo-spectral scheme for the GR-KdV equation in 2D. In addition, to the best of our knowledge, there has been no reference considering an error estimate for  the Fourier spectral schemes of the GR-KdV equation. Thus, another purpose of this paper is to establish an a priori estimate for the proposed Fourier pseudo-spectral
scheme in discrete $L^{\infty}$-norm.

The outline of this paper is
organized as follows. In Section 2, a new semi-norm equivalence between the Fourier pseudo-spectral method and the finite difference method is first established. A semi-discrete system of the GR-KdV equation \eqref{GR-KDV_eq:1.1}, which inherits the semi-discrete momentum conservation law, is then presented by using the standard Fourier pseudo-spectral method in space. In Section 3, a fully discrete scheme is obtained and we show that the resulting scheme is momentum-preserving and uniquely solvable. An a priori estimate is established for the proposed
scheme in discrete $L^{\infty}$-norm in Section 4. Some numerical experiments are presented in Section 5. We draw some conclusions in Section 6.

\section{Structure-preserving spatial discretization}
Let $\Omega_{h}=\{(x_{j_1},y_{j_2})|x_{j_{1}}=j_{1}h_{1},y_{j_{2}}=j_{2}h_{2};\ 0\leq j_{r}\leq N_{r}, r=1,2\}$ be a partition of $\Omega$
with mesh sizes $h_1=\frac{l_{1}}{N_{1}}$ and $h_2=\frac{l_{2}}{N_{2}}$, respectively, where $N_{1}$ and $N_{2}$ are two even numbers. Denote \begin{align*}
J_{h}=\{(j_1,j_2)|0\leq j_{r}\leq N_{r}, r=1,2\},\ J_{h}^{'}=\{(j_1,j_2)|0\leq j_{r}\leq N_{r}-1, r=1,2\}.
\end{align*}
A discrete mesh function ${ U}_{j_1,j_2},\ (j_1,j_2)\in \mathbb{Z}\times\mathbb{Z}$ is said to satisfy periodic boundary conditions if and only if
\begin{align}\label{GR-KDV-PBS}
&x-\text{periodic}:\  U_{j_1,j_2}=U_{j_1+N_1,j_2}\ \text{and}\ y-\text{periodic}:\ U_{j_1,j_2}=U_{j_1,j_2+N_2}.
\end{align}

Let
\begin{align*}
\mathbb{ V}_{h}:&=\big\{{\bm U}|{\bm U}=(U_{0,0},U_{1,0},\cdots,U_{N_{1}-1,0},U_{0,1},U_{1,1},\cdots, U_{N_{1}-1,1},\cdots,U_{0,N_{2}-1},U_{1,N_{2}-1},\\
 &\cdots,U_{N_1-1,N_{2}-1})^{T}\big\}
\end{align*}
be the space of mesh functions defined on $\Omega_h$ and satisfy the periodic boundary conditions \eqref{GR-KDV-PBS}.
Subsequently, the discrete difference operators and norms will be defined in an appropriate way. We first introduce some discrete difference operators for any mesh function ${\bm U}\in \mathbb{V}_h$, as follows:
\begin{align*}
&\delta_{x}^{+} U_{j_1,j_2}=\frac{U_{j_{1}+1,j_{2}}-U_{j_{1},j_{2}}}{h_{1}},\ \delta_{y}^{+} U_{j_1,j_2}=\frac{U_{j_{1},j_{2}+1}-U_{j_{1},j_{2}}}{h_{2}},\\
& \delta_{x}^{-} U_{j_1,j_2}=\frac{U_{j_{1},j_{2}}-U_{j_{1}-1,j_{2}}}{h_{1}},\ \delta_{y}^{-} U_{j_1,j_2}=\frac{U_{j_{1},j_{2}}-U_{j_{1},j_{2}-1}}{h_{2}},\\
%& \widehat{A}_{t}U_{\vec{j}}^{n}=\frac{U_{\vec{j}}^{n+1}+U_{\vec{j}}^{n-1}}{2},\ U_{\vec{j}}^{n+\frac{1}{2}}=\frac{U_{\vec{j}}^{n+1}+U_{\vec{j}}^{n}}{2},\
%\hat{U}_{\vec{j}}^{n+\frac{1}{2}}=\frac{3U_{\vec{j}}^{n}-U_{\vec{j}}^{n-1}}{2},\
&\nabla_h U_{j_1,j_2}=\Big(\delta_{x}^+U_{j_1,j_2},\ \delta_{y}^+ U_{j_1,j_2}\Big)^T,\ \Delta_h U_{j_1,j_2}=(\delta_{x}^{+}\delta_{x}^{-}+\delta_{y}^{+}\delta_{y}^{-})U_{j_1,j_2},\ (j_1,j_2)\in J_{h}^{'}.
\end{align*}
%\begin{align*}
%&\langle {\bm U},{\bm V}\rangle_{h}=h_1h_2\sum_{j_1=0}^{N_1-1}\sum_{j_2=0}^{N_2-1}U_{j_1,j_2}V_{j_1,j_2},\ \forall\  {\bm U},{\bm V}\in\mathbb{ V}_{h}.
%%h_{\Delta}\sum_{j_1=0}^{N_1-1}\sum_{j_2=0}^{N_2-1}\sum_{j_3=0}^{N_3-1}U_{j_1,j_2,j_3}\overline{V}_{j_1,j_2,j_3}=
%\end{align*}
Then, for any ${\bm U}$ and ${\bm V}$ in $\mathbb{V}_h$, we define the discrete inner product and notions as
%The discrete $L^{2}$-norm of ${\bm U}\in\mathbb{V}_h$ and its difference quotients are defined, respectively, as
\begin{align*}
&\langle {\bm U},{\bm V}\rangle_{h}=h_1h_2\sum_{j_1=0}^{N_1-1}\sum_{j_2=0}^{N_2-1}U_{j_1,j_2}V_{j_1,j_2},\ ||{\bm U}||_{h}^2=h_1h_2\sum_{j_1=0}^{N_1-1}\sum_{j_2=0}^{N_2-1}|U_{j_1,j_2}|^2,\\
 &||{\delta_{x}^+\bm U}||_h^2=h_1h_2\sum_{j_1=0}^{N_1-1}\sum_{j_2=0}^{N_2-1}|\delta_{x}^+U_{j_1,j_2}|^2,\ ||{\delta_{y}^+\bm U}||_h^2=h_1h_2\sum_{j_1=0}^{N_1-1}\sum_{j_2=0}^{N_2-1}|\delta_{y}^+U_{j_1,j_2}|^2,\\
 &||\nabla_h{\bm U}||_h^2=||{\delta_{x}^+\bm U}||_h^2+||{\delta_{y}^+\bm U}||_h^2,\ ||\Delta_h{\bm U}||_h^2=h_1h_2\sum_{j_1=0}^{N_1-1}\sum_{j_2=0}^{N_2-1}|\Delta_hU_{j_1,j_2}|^2.
\end{align*}
We also define discrete $H_h^2$ and $L^{\infty}$-norms as
\begin{align*}
& ||{\bm U}||_{H_h^2}^2={||{\bm U}||_h^2+||\nabla_h{\bm U}||_h^2+||\Delta_h{\bm U}||_h^2},\ ||{\bm U}||_{h,\infty}=\max\limits_{(j_1,j_2)\in J_{h}^{'}}|U_{j_1,j_2}|.%,\ \forall\ {\bm U}\in\mathbb{V}_h.
\end{align*}
We note that the discrete norms $||\delta_x^+{\bm U}||_h$, $||\delta_y^+{\bm U}||_h$, $||\nabla_h{\bm U}||_h$ and $||\Delta_h{\bm U}||_h$ defined above are semi-norms.
In addition, we denote $`\cdot$' as the element product of vectors ${\bm U}, {\bm V}\in\mathbb{V}_h$, that is,
\begin{align*}
%&{\bm U}=({\hat {\bm U}}_{0,0}^{T},\cdots,{\hat{{\bm U}}}_{N_{2}-1,0}^{T},{\hat {\bm U}}_{0,1}^{T},\cdots,{\hat {\bm U}}_{N_{2}-1,1}^{T},\cdots,{\hat {\bm U}}_{0,N_{3}-1}^{T},\cdots,{\hat {\bm U}}_{N_{2}-1,N_{3}-1}^{T})^{T},\\
%&{\hat {\bm U}}_{j_2,j_3}=(U_{0,j_2,j_3},\cdots,U_{N_{1}-1,j_2,j_3})^{T}.
{\bm U}\cdot {\bm V}=&\big(U_{0,0}V_{0,0},\cdots,U_{N_1-1,0}V_{N_1-1,0},\cdots,U_{0,N_2-1}V_{0,N_2-1},\cdots,U_{N_1-1,N_2-1}V_{N_1-1,N_2-1}\big)^{T}.
\end{align*}
For brevity, we denote $\underbrace{{\bm U}\cdot...\cdot {\bm U}}_{k}$ as ${\bm U}^k$.
\begin{den}\label{GR-KDV:den2.1} {In this paper, for any matrices ${\bm A}=(a_{j,k})_{pq}$ and ${\bm B}=(b_{j,k})_{rm}$, where $p,\ q,\ r,\ m$ are nonnegative integers, the Kronecker product ${\bm A}\otimes{\bm B}$ is a $pr\times qm$ block matrix defined by $$
{\bm A}\otimes{\bm B}=\left(\begin{array}{lllll}
a_{1,1}{\bm B}& a_{1,2}{\bm B}  &\cdots &a_{1,q}{\bm B} \\
              a_{2,1}{\bm B}& a_{2,2}{\bm B}  &\cdots &a_{2,q}{\bm B} \\
              \vdots &\vdots  &\ddots &\vdots\\
              a_{p,1}{\bm B}&  a_{p,2}{\bm B}& \cdots & a_{p,q}{\bm B}\\
\end{array}
\right).
$$ }
\end{den}
\begin{coro}\label{GR-KDV:cor:2.1} {According to the definition \ref{GR-KDV:den2.1}, we can show that, for any matrices ${\bm A}=(a_{j,k})_{pm}$, ${\bm B}=(b_{j,k})_{ml}$, ${\bm C}=(c_{j,k})_{rq}$ and ${\bm D}=(d_{j,k})_{qs}$, where $p,\ q,\ r,\ m,\ l,\ s$ are nonnegative integers, the Kronecker product $\otimes$ satisfies
\begin{align*}
({\bm A}\otimes {\bm C})({\bm B}\otimes {\bm D})={\bm A}{\bm B}\otimes{\bm C}{\bm D}.
\end{align*}}

\end{coro}

\begin{rmk} { According to the definition \ref{GR-KDV:den2.1} and Corollary \ref{GR-KDV:cor:2.1},  $||\nabla_h{\bm U}||_h^2$ and $||\Delta_h{\bm U}||_h^2$ can be rewritten as
\begin{align*}
&||\nabla_h{\bm U}||_h^2={\langle-\big({\bm I}_{N_{2}}\otimes{\bm B}_{1}+{\bm B}_{2}\otimes {\bm I}_{N_{1}}\big){\bm U},{\bm U}\rangle_{h}}:={\langle-{\bm \Delta}_h{\bm U},{\bm U}\rangle_{h}},\\
&||\Delta_h{\bm U}||_h^2={\langle\big({\bm I}_{N_{2}}\otimes{\bm B}_{1}^2+2{\bm B}_{2}\otimes{\bm B}_1+{\bm B}_{2}^2\otimes {\bm I}_{N_{1}}\big){\bm U},{\bm U}\rangle_h}:={\langle{\bm \Delta}^2_h{\bm U},{\bm U}\rangle_{h}},
\end{align*}
where $I_{N_r}$ is an $N_r\times N_r$ identity matrix, and
\begin{align*}
{\bm B}_r=\frac{1}{h_r^2}\left(\begin{array}{ccccc}
              -2& 1 &0 &\cdots &1 \\
              1& -2 &1& \cdots &0\\
              \vdots &\vdots &\ddots& \ddots &\vdots\\
              0&0& \cdots&-2&1\\
              1&  0& \cdots&1 & -2\\
             \end{array}
\right)_{N_r\times N_r},r=1,2.
\end{align*}
Here, ${\bm B}_r, r=1,2$ is the usual finite diffidence discretization of the second derivative, by taking
into account of the periodic boundary conditions. }
\end{rmk}

%By simple calculation, we find that
%\begin{align*}
%|{\bm u}|_{1,h}=\langle{\bm B}{\bm u}, {\bm u}\rangle_{J}^{\frac{1}{2}},
%\end{align*}
%with
%\begin{align*}
%{\bm B}=\frac{1}{h^2}\left(\begin{array}{ccccc}
%              2& -1 &0 &\cdots &-1 \\
%              -1& 2 &-1& \cdots &0\\
%              \vdots &\vdots &\ddots& \ddots &\vdots\\
%              0&0& \cdots&2&-1\\
%              -1&  0& \cdots&-1 & 2\\
%             \end{array}
%\right)_{J\times J}.
%\end{align*}
%%We note that $|{\bm u}|_{1,h}$ is a semi-norm.
%\begin{rmk}
%
%\end{rmk}
\subsection{Fourier pseudo-spectral method and some useful lemmas}
Let
\begin{align*}
&S_{N}^{''}=\text{span}\{g_{j_1}(x)g_{j_2}(y),\ 0\leq j_r\leq N_r-1, r=1,2\},
\end{align*}
be the interpolation space, where $g_{j_1}(x)$ and $g_{j_2}(y)$ are trigonometric polynomials of degree $N_{1}/2$ and $N_{2}/2$,
given, respectively, by
\begin{align*}
  &g_{j_1}(x)=\frac{1}{N_{1}}\sum_{l=-N_{1}/2}^{N_{1}/2}\frac{1}{a_{l}}e^{\text{i}l\mu_{1} (x-x_{j_1})},\ g_{j_2}(y)=\frac{1}{N_{2}}\sum_{q=-N_{2}/2}^{N_{2}/2}\frac{1}{b_{q}}e^{\text{i}q\mu_{2} (y-y_{j_2})},
\end{align*}
with $a_{l}=\left \{
 \aligned
 &1,\ |l|<\frac{N_1}{2},\\
 &2,\ |l|=\frac{N_1}{2},
 \endaligned
 \right.,\ b_{q}=\left \{
 \aligned
 &1,\ |q|<\frac{N_2}{2},\\
 &2,\ |q|=\frac{N_2}{2},
 \endaligned
 \right.$, $\mu_{r}=\frac{2\pi}{l_{r}}$ and $\ 0\le j_r\le N_r-1,\ r=1,2$.
We define the interpolation operator $I_{N}: C(\Omega)\to S_{N}^{''}$ as \cite{CQ01}:
\begin{align*}
I_{N}U(x,y,t)=\sum_{j_1=0}^{N_{1}-1}\sum_{j_2=0}^{N_{2}-1}U_{j_1,j_2}(t)g_{j_1}(x)g_{j_2}(y),
\end{align*}
where $U_{j_1,j_2}(t)=U(x_{j_1},y_{j_2},t)$. %and its vector form is denoted by
%\begin{align*}
%{\bm u}=&(u_{0,0,0},u_{1,0,0},\cdots,u_{N_{1}-1,0,0},u_{0,1,0},u_{1,1,0},\cdots, u_{N_{1}-1,1,0},\cdots,u_{0,N_{2}-1,0},u_{1,N_{2}-1,0},\cdots,\\
% &u_{N_1-1,N_{2}-1,0},u_{0,0,1},u_{1,0,1}\cdots,u_{N_1-1,0,1},\cdots,u_{0,N_2-1,N_3-1},\cdots,u_{N_1-1,N_2-1,N_3-1})^{T}.
%\end{align*}

Taking the derivative with respect to $x$, and then evaluating the resulting expressions at the collocation points ($x_{j_1},y_{j_2}$), where $0\le j_r\le N_r-1,\ r=1,2$, we have
\begin{align}\label{GR-KDV-2.2}
\frac{\partial^{s_{1}} I_{N}U(x_{j_1},y_{j_2},t)}{\partial x^{s_1}}
&=\sum_{j=0}^{N_{1}-1}U_{j,j_2}(t)\frac{d^{s_1}g_{j}(x_{j_1})}{dx^{s_1}}=[({\bm I}_{N_{2}}\otimes{\bm D}_{s_1}^{x}){\bm U}]_{j_1,j_2}, {\bm U}\in\mathbb{V}_h,%{{N_{1}N_{2}(j_3-1)+N_{1}(j_2-1)+j_1}}
\end{align}
where ${\bm D}_{s_1}^{x}$ is an $N_{1}\times N_1$ matrix, with elements given by
\begin{align*}
({\bm D}_{s_1}^{x})_{j_1,j}=\frac{d^{s_1}g_{j}(x_{j_1})}{dx^{s_1}},\ 0\le j_1,j\le N_1-1,
\end{align*}
and $[({\bm I}_{N_{2}}\otimes{\bm D}_{s_1}^{x}){\bm U}]_{j_1,j_2}$ represents the $(N_{1}j_2+j_1+1)$-th component of the vector $({\bm I}_{N_{2}}\otimes{\bm D}_{s_1}^{x}){\bm U}$. For brevity, the notation is still be adopted in subsequent sections. Similarly, we can obtain
\begin{align}\label{GR-KDV-2.3}
%&:=[{\bm D}_{1}{\bm U}]_{N_{x}N_{y}(m-1)+N_{x}(k-1)+j},\\
\frac{\partial^{s_2} I_{N}U(x_{j_1},y_{j_2},t)}{\partial y^{s_2}}&=\sum_{k=0}^{N_{2}-1}U_{j_1,k}\frac{d^{s_2}g_{k}(y_{j_2})}{dy^{s_2}}
=[({\bm D}_{s_2}^{y}\otimes {\bm I}_{N_{1}}){\bm U}]_{j_1,j_2},\ {\bm U}\in\mathbb{V}_h,
%&:=[{\bm D}_{2}{\bm U}]_{{N_{x}N_{y}(m-1)+N_{x}(k-1)+j}},\\
\end{align}
where %$\otimes$ is Kronecker product, $I_{N_{r}},\ r=1,2,3$ is the identity matrix of dimension $N_{r}\times N_{r}$ and
 ${\bm D}_{s_2}^{y}$ is an $N_{2}\times N_2$ matrix with elements given by
\begin{align*}
 ({\bm D}_{s_2}^{y})_{j_2,k}=\frac{d^{s_2}g_{k}(y_{j_2})}{dy^{s_2}},\ 0\le j_2,k\le N_2-1.
\end{align*}
In particular, for first and second derivatives, we have, respectively,
\begin{align*}
\frac{\partial I_{N}U(x_{j_1},y_{j_2},t)}{\partial x}
=[({\bm I}_{N_{2}}\otimes{\bm D}_{1}^{x}){\bm U}]_{j_1,j_2},\
\frac{\partial I_{N}U(x_{j_1},y_{j_2},t)}{\partial y}
=[({\bm D}_{1}^{y}\otimes {\bm I}_{N_{1}}){\bm U}]_{j_1,j_2},
\end{align*}
and
\begin{align*}
\frac{\partial^{2} I_{N}U(x_{j_1},y_{j_2},t)}{\partial x^{2}}
=[({\bm I}_{N_{2}}\otimes{\bm D}_{2}^{x}){\bm U}]_{j_1,j_2},\
\frac{\partial^{2} I_{N}U(x_{j_1},y_{j_2},t)}{\partial y^{2}}
=[({\bm D}_{2}^{y}\otimes {\bm I}_{N_{1}}){\bm U}]_{j_1,j_2},
\end{align*}
where ${\bm D}_{1}^{x}$ and ${\bm D}_{1}^{y}$ are real skew-symmetric matrices, and ${\bm D}_{2}^{x}$ and ${\bm D}_{2}^{y}$ are real symmetric matrices, respectively.
\begin{rmk} { With noting definition  \ref{GR-KDV:den2.1} and $$
({\bm I}_{N_{2}}\otimes{\bm D}_{s_1}^{x}){\bm U}=vec({\bm D}_{s_1}^{x}\mathcal{U}),\ ({\bm D}_{s_2}^{y}\otimes {\bm I}_{N_{1}}){\bm U}=vec(\mathcal{U}({\bm D}_{s_2}^{y})^T),$$
where $\mathcal{U}$ is an $N_1$-by-$N_2$ matrix whose elements
    are taken columnwise from ${\bm U}$ and the vec operator stacks the columns of a matrix one underneath the other to form a single vector.
    Then, it is clear to see that
    \begin{align*}
    [({\bm I}_{N_{2}}\otimes{\bm D}_{s_1}^{x}){\bm U}]_{j_1,j_2}=[vec({\bm D}_{s_1}^{x}\mathcal{U})]_{j_1,j_2}=\sum_{k=0}^{N_{1}-1}({\bm D}_{s_1}^{x})_{j_1,k}U_{k,j_2},
    \end{align*}
    and
    \begin{align*}
    [({\bm D}_{s_2}^{y}\otimes {\bm I}_{N_{1}}){\bm U}]_{j_1,j_2}=[vec(\mathcal{U}({\bm D}_{s_2}^{y})^T)]_{j_1,j_2}=\sum_{k=0}^{N_{2}-1}U_{j_1,k}({\bm D}_{s_2}^{y})_{j_2,k}.
    \end{align*}}
\end{rmk}
Now, for any mesh function ${\bm U}\in\mathbb{V}_{h}$, we define three new semi-norms induced by the spectral differential matrices as
\begin{align*}
 &|{\bm U}|_{1,h}^2={\langle-({\bm I}_{N_{2}}\otimes({\bm D}_{1}^{x})^2+ ({\bm D}_{1}^{y})^2\otimes {\bm I}_{N_{1}}){\bm U},{\bm U}\rangle_{h}},\ ||\mathbb{L}_h{\bm U}||_{h}^2={\langle\mathbb{L}_h{\bm U}, \mathbb{L}_h{\bm U}\rangle},\\
 &|{\bm U}|_{2,h}^2={\langle({\bm I}_{N_{2}}\otimes{\bm D}_{2}^{x}+{\bm D}_{2}^{y}\otimes {\bm I}_{N_{1}}){\bm U},({\bm I}_{N_{2}}\otimes{\bm D}_{2}^{x}+{\bm D}_{2}^{y}\otimes {\bm I}_{N_{1}}){\bm U}\rangle_{h}},
\end{align*}
where $\mathbb{L}_h={\bm I}_{N_{2}}\otimes{\bm D}_{1}^{x}+ {\bm D}_{1}^{y}\otimes {\bm I}_{N_{1}}.$
\begin{lem}\rm\label{GR-KDV-lem2.1} For any mesh function ${\bm U}\in\mathbb{V}_{h}$, we have
\begin{align*}
||\mathbb{L}_h{\bm U}||_{h}\le \sqrt{2}|{\bm U}|_{1,h}.
\end{align*}
\begin{prf}\rm { According to Corollary \ref{GR-KDV:cor:2.1}, we have }
\begin{align}\label{GR-KDV:eq:2.4}
{||\mathbb{L}_h{\bm U}||_{h}^2}&{=\langle -\mathbb{L}_h^2{\bm U},{\bm U}\rangle_h}\nonumber\\
&{=\langle-({\bm I}_{N_{2}}\otimes({\bm D}_{1}^{x})^2+2{\bm D}_{1}^{y}\otimes {\bm D}_{1}^{x}+({\bm D}_{1}^{y})^2\otimes I_{N_1}){\bm U},{\bm U}\rangle_{h}}\nonumber\\
%&\textcolor{blue}{=|{\bm U}|_{1,h}^2-2\langle {\bm D}_{1}^{y}\otimes {\bm D}_{1}^{x}{\bm U},{\bm U}\rangle_{h}}\\
&{=|{\bm U}|_{1,h}^2-2\langle({\bm I}_{N_2}\otimes {\bm D}_{1}^{x}){\bm U},({\bm D}_{1}^{y}\otimes {\bm I}_{N_1})^{T}{\bm U}\rangle_{h}}\nonumber\\
&{=|{\bm U}|_{1,h}^2+2\langle({\bm I}_{N_2}\otimes {\bm D}_{1}^{x}){\bm U},({\bm D}_{1}^{y}\otimes {\bm I}_{N_1}){\bm U}\rangle_{h}}.
\end{align}
{ It follows from the Cauchy-Schwarz inequality and mean value inequality  that}
\begin{align}\label{GR-KDV:eq:2.5}
{2\langle({\bm I}_{N_2}\otimes {\bm D}_{1}^{x}){\bm U},({\bm D}_{1}^{y}\otimes {\bm I}_{N_1}){\bm U}\rangle_{h}}&{\le 2||({\bm I}_{N_2}\otimes {\bm D}_{1}^{x}){\bm U}||_h||({\bm D}_{1}^{y}\otimes {\bm I}_{N_1}){\bm U}||_h}\nonumber\\
&{\le||({\bm I}_{N_2}\otimes {\bm D}_{1}^{x}){\bm U}||_h^2+||({\bm D}_{1}^{y}\otimes {\bm I}_{N_1}){\bm U}||_h^2}\nonumber\\
&{=|{\bm U}|_{1,h}^2,}
\end{align}
Combining \eqref{GR-KDV:eq:2.4} and \eqref{GR-KDV:eq:2.5}, we have
\begin{align*}
||\mathbb{L}_h{\bm U}||_{h}\le \sqrt{2}|{\bm U}|_{1,h}.
\end{align*}
This completes the proof.\qed
\end{prf}
\end{lem}
\begin{lem}\rm \label{GR-KDV-lem2.2}\cite{GCW14,HNO06} For the matrices ${\bm B}_{r},r=1,2$ and ${\bm D}_{2}^{w},w=x,y$, the following results hold
\begin{align*}
&{\bm B}_{r}={\bm F}_{N_r}^{H}\Lambda_{r}{\bm F}_{N_r},\\
&{\bm D}_{1}^{x}={\bm F}_{N_1}^{H}\Lambda_{3}{\bm F}_{N_1},\\
&{\bm D}_{1}^{y}={\bm F}_{N_2}^{H}\Lambda_{4}{\bm F}_{N_2},\\
&{\bm D}_{2}^{x}={\bm F}_{N_1}^{H}\Lambda_{5}{\bm F}_{N_1},\\
&{\bm D}_{2}^{y}={\bm F}_{N_2}^{H}\Lambda_{6}{\bm F}_{N_2},
\end{align*}
where  ${\bm F}_{N_r},\ r=1,2$, is the discrete Fourier transform matrix with elements
$\big({\bm F}_{N_r}\big)_{j,k}=\frac{1}{\sqrt{N_r}}e^{-\text{\rm i}jk\frac{2\pi}{N_r}},$ ${\bm F}_{N_r}^{H}$ is the conjugate transpose matrix of ${\bm F}_{N_r}$ and
\begin{align*}
&\Lambda_{r}=\text{\rm diag}\Big[\lambda_{B_{r},0},\lambda_{B_{r},1},\cdots,\lambda_{B_{r},N_{r}-1}\Big],\ \lambda_{B_{r},j}=-\frac{4}{h_{r}^{2}}\sin^{2}\frac{j\pi}{N_{r}},\\
&\Lambda_{3}=\text{\rm diag}\Big[\lambda_{D_{1}^{x},0},\lambda_{D_{1}^{x},1},\cdots,\lambda_{D_{1}^{x},N_{1}-1}\Big],\ \lambda_{D_{1}^{x},j}=
 \left \{
 \aligned
 &\text{i}j\mu_{1},\ &0\leq j\leq N_{1}/2-1,\\
 &0,\ &j=N_1/2,\\
 &\text{i}(j-N_1)\mu_1,\ &N_1/2<j<N_1,
 \endaligned
 \right.\\
 &\Lambda_{4}=\text{\rm diag}\Big[\lambda_{D_{1}^{y},0},\lambda_{D_{1}^{y},1},\cdots,\lambda_{D_{1}^{y},N_{2}-1}\Big],\ \lambda_{D_{1}^{y},j}=
 \left \{
 \aligned
 &\text{i}j\mu_{2},\ &0\leq j\leq N_{2}/2-1,\\
 &0,\ &j=N_2/2,\\
 &\text{i}(j-N_2)\mu_2,\ &N_2/2<j<N_2,
 \endaligned
 \right.\\
&\Lambda_{5}=\text{\rm diag}\Big[\lambda_{D_{2}^{x},0},\lambda_{D_{2}^{x},1},\cdots,\lambda_{D_{2}^{x},N_{1}-1}\Big],\ \lambda_{D_{2}^{x},j}=
 \left \{
 \aligned
 &-(j\mu_{1})^{2},\ &0\leq j\leq N_{1}/2,\\
 &-\big((j-N_1)\mu_1\big)^{2},\ &N_1/2<j<N_1,
 \endaligned
 \right.\\
 &\Lambda_{6}=\text{\rm diag}\Big[\lambda_{D_{2}^{y},0},\lambda_{D_{2}^{y},1},\cdots,\lambda_{D_{2}^{y},N_{2}-1}\Big],\ \lambda_{D_{2}^{y},j}=
 \left \{
 \aligned
 &-(j\mu_{2})^{2},\ &0\leq j\leq N_{2}/2,\\
 &-\big((j-N_2)\mu_2\big)^{2},\ &N_2/2<j<N_2.
 \endaligned
 \right.
\end{align*}
In addition, the following inequalities hold \cite{GWWC17}
\begin{align}\label{GR-KDV:2.1}
%&0\leq-\frac{4}{\pi^{2}}\lambda_{(D_{1}^{x})^2,j}\leq -\lambda_{B_{1},j},\ 0\leq j\leq N_1-1,\\\label{GR-KDV:2.01}
%&0\leq-\frac{4}{\pi^{2}}\lambda_{(D_{1}^{y})^2,j}\leq -\lambda_{B_{2},j},\ 0\leq j\leq N_2-1,\\\label{GR-KDV:2.1}
&0\leq-\frac{4}{\pi^{2}}\lambda_{(D_{1}^{x})^2,j}\leq-\frac{4}{\pi^{2}}\lambda_{D_{2}^{x},j}\leq -\lambda_{B_{1},j}\leq -\lambda_{D_{2}^{x},j},\ 0\leq j\leq N_1-1,\\\label{GR-KDV:2.2}
&0\leq-\frac{4}{\pi^{2}}\lambda_{(D_{1}^{y})^2,j}\leq-\frac{4}{\pi^{2}}\lambda_{D_{2}^{y},j}\leq -\lambda_{B_{2},j}\leq -\lambda_{D_{2}^{y},j},\ 0\leq j\leq N_2-1,\\\label{GR-KDV:2.3}
&0\leq\frac{16}{\pi^{4}}\lambda_{D_{2}^{x},j}^2\leq \lambda_{B_{1},j}^2\leq \lambda_{D_{2}^{x},j}^2,\ 0\leq j\leq N_1-1,\\\label{GR-KDV:2.4}
&0\leq\frac{16}{\pi^{4}}\lambda_{D_{2}^{y},j}^2\leq \lambda_{B_{2},j}^2\leq \lambda_{D_{2}^{y},j}^2,\ 0\leq j\leq N_2-1.
\end{align}
\end{lem}
%Now, we give the following equivalence between $|{\bm U}|_{1,h}$ and $|{\bm U}|_{h}$.
\begin{lem}\rm \label{GR-KDV-lem2.3} For any mesh function ${\bm U} \in \mathbb{V}_{h}$, we have
\begin{align*}%\label{eq:2.11}
&||\Delta_h{\bm U}||_{h}\leq |{\bm U}|_{2,h}\leq\frac{\pi^2}{4}||\Delta_h{\bm U}||_{h}.
\end{align*}
\end{lem}
\begin{prf}\rm
We denote
\begin{align*}
I^{2}:&=|{\bm U}|_{2,h}^{2}\\
&=\langle\big({\bm I}_{N_{2}}\otimes({\bm D}_{2}^{x})^2\big){\bm U},{\bm U}\rangle_{h}+2\langle\big({\bm D}_{2}^{y}\otimes {\bm D}_{2}^{x}\big){\bm U},{\bm U}\rangle_{h}+\langle\big(({\bm D}_{2}^{y})^2\otimes {\bm I}_{N_{1}}\big){\bm U},{\bm U}\rangle_{h}\nonumber\\
&:=I_{1}^{2}+2I_{2}^{2}+I_{3}^{2},
\end{align*}
and
\begin{align*}
J^{2}:&=||\Delta_{h}{\bm U}||_h^2=\langle\Delta_{h}{\bm U},\Delta_{h}{\bm U}\rangle_{h}\nonumber\\
&=\langle\big({\bm I}_{N_{2}}\otimes{\bm B}_{1}^2\big){\bm U},{\bm U}\rangle_{h}+2\langle\big({\bm B}_{2}\otimes {\bm B}_{1}\big){\bm U},{\bm U}\rangle_{h}+\langle\big({\bm B}_{2}^2\otimes {\bm I}_{N_{1}}\big){\bm U},{\bm U}\rangle_{h}\\
&:=J_1^{2}+2J_2^{2}+J_3^{2}.
\end{align*}
With Lemma \ref{GR-KDV-lem2.2}, we can obtain
\begin{align*}
I_{1}^{2}&=\langle\big({\bm I}_{N_{2}}\otimes{\bm F}_{N_1}^{H}\Lambda_{5}^2{\bm F}_{N_1}\big){\bm U},{\bm U}\rangle_{h}
=\langle\big({\bm I}_{N_{2}}\otimes\Lambda_{5}^2\big)\widetilde{{\bm U}},\widetilde{{\bm U}}\rangle_{h}\nonumber\\
&=h_1h_2\sum_{j_1=0}^{N_1-1}\sum_{j_2=0}^{N_2-1}\big(\lambda_{D_{2}^{x},j_1}\big)^2|\widetilde{{U}}_{j_1,j_2}|^{2},
\end{align*}
\begin{align*}
I_{2}^{2}&=\langle\big( {\bm F}_{N_2}^{H}\Lambda_{6}{\bm F}_{N_2}\otimes {\bm F}_{N_1}^{H}\Lambda_{5}{\bm F}_{N_1}\big){\bm U},{\bm U}\rangle_{h}
=\langle\big(\Lambda_{6}\otimes\Lambda_{5}\big)\widetilde{{\bm U}},\widetilde{{\bm U}}\rangle_{h}\nonumber\\
&=h_1h_2\sum_{j_1=0}^{N_1-1}\sum_{j_2=0}^{N_2-1}\big(\lambda_{D_{2}^{x},j_1}\lambda_{D_{2}^{y},j_2}\big)|\widetilde{{U}}_{j_1,j_2}|^{2},
\end{align*}
and
\begin{align*}
I_{3}^{2}&=\langle\big( {\bm F}_{N_2}^{H}\Lambda_{6}^2{\bm F}_{N_2}\otimes {\bm I}_{N_{1}}\big){\bm U},{\bm U}\rangle_{h}
=\langle\big(\Lambda_{6}^2\otimes{\bm I}_{N_{1}}\big)\widetilde{{\bm U}},\widetilde{{\bm U}}\rangle_{h}\nonumber\\
&=h_1h_2\sum_{j_1=0}^{N_1-1}\sum_{j_2=0}^{N_2-1}\big(\lambda_{D_{2}^{y},j_2}\big)^2|\widetilde{{U}}_{j_1,j_2}|^{2},
\end{align*}
where $\widetilde{{\bm U}}=\big({\bm F}_{N_2}\otimes{\bm F}_{N_1}\big){\bm U}$ and $\widetilde{{U}}_{j_1,j_2}=({\bm F}_{N_1}\mathcal{U}{\bm F}_{N_2}^T)_{j_1,j_2}$.
Similarly, we can deduce
\begin{align*}
&J_1^{2}=h_1h_2\sum_{j_1=0}^{N_1-1}\sum_{j_2=0}^{N_2-1}\big(\lambda_{B_{1},j_1}\big)^2|\widetilde{{ U}}_{j_1,j_2}|^{2},\ J_2^{2}=h_1h_2\sum_{j_1=0}^{N_1-1}\sum_{j_2=0}^{N_2-1}\big(\lambda_{B_{1},j_1}\lambda_{B_{2},j_2}\big)|\widetilde{{ U}}_{j_1,j_2}|^{2},\\
&J_3^{2}=h_1h_2\sum_{j_1=0}^{N_1-1}\sum_{j_2=0}^{N_2-1}\big(\lambda_{B_{2},j_2}\big)^2|\widetilde{{ U}}_{j_1,j_2}|^{2}.
\end{align*}
%\begin{align*}
%J_2^{2}=h_{\Delta}\sum_{j_1=0}^{N_1-1}\sum_{j_2=0}^{N_2-1}\sum_{j_3=0}^{N_3-1}\big(-\lambda_{B_{2},j_2}\big)|\widetilde{{U}}_{j_1,j_2,j_3}|^{2},
%\end{align*}
%and
%\begin{align*}
%J_3^{2}=h_{\Delta}\sum_{j_1=0}^{N_1-1}\sum_{j_2=0}^{N_2-1}\sum_{j_3=0}^{N_3-1}\big(-\lambda_{B_{3},j_3}\big)|\widetilde{{ U}}_{j_1,j_2,j_3}|^{2}.
%\end{align*}
With the use of \eqref{GR-KDV:2.1}-\eqref{GR-KDV:2.4}, we have
\begin{align*}%\label{3d_NLS:2.23}
\frac{16}{\pi^4}I_{r}^{2}\leq J_r^{2}\leq I_r^{2},\  r=1,2,3,
\end{align*}
which implies that
\begin{align}\label{GR-KDV:2.5}
J_{r}^{2}\leq I_r^{2}\leq \frac{\pi^4}{16}J_r^{2}, \ r=1,2,3.
\end{align}
With \eqref{GR-KDV:2.5}, we can get
\begin{align*}%\label{3d_NLS:2.25}
J^{2}\leq I^{2}\leq \frac{\pi^4}{16}J^{2},
\end{align*}
that is,
\begin{align*}%\label{3d_NLS:2.26}
||\Delta_h{\bm U}||_{h}\leq |{\bm U}|_{2,h}\leq\frac{\pi^2}{4}||\Delta_h{\bm U}||_{h}.
\end{align*}
This completes the proof.
\qed
\end{prf}
\begin{lem}\rm \label{GR-KDV-lem2.4} For any mesh function ${\bm U} \in \mathbb{V}_{h}$, we have
\begin{align*}%\label{eq:2.11}
&|{\bm U}|_{1,h}\leq\frac{\pi}{2}||\nabla_h{\bm U}||_{h}.
\end{align*}
\end{lem}
The proof is similar to the Lemma \ref{GR-KDV-lem2.3}. For brevity, we omit it.
\begin{lem}\rm \cite{LS10,QS12,zhou90} \label{GR-KDV-lem2.5} For any mesh function ${\bm U}\in \mathbb{V}_h$, we have
\begin{align*}
&||\nabla_h {\bm U}||_h^2\le ||{\bm U}||_h||\Delta_h {\bm U}||_h,\\
&||{\bm U}||_{h,\infty}^2\le C||{\bm U}||_h\big(||\Delta_h {\bm U}||_h+||{\bm U}||_h\big).
\end{align*}

\end{lem}
\subsection{Momentum-preserving spatial semi-discretization}
Eq. \eqref{GR-KDV_eq:1.1} can be rewritten as the following equivalent form
\begin{align}\label{GR-KDV_eq:2.6}
 u_t+\Delta^2u_t+\Delta u_x+\mathbb{L}u+\frac{1}{p+2}(u^p\mathbb{L} +\mathbb{L}u^p)u=0,
\end{align}
where $(u^p\mathbb{L} +\mathbb{L}u^p)u=u^p\mathbb{L}u +\mathbb{L}(u^pu)$, which is applicable for the discrete version.

Applying the Fourier pseudo-spectral method to the system \eqref{GR-KDV_eq:2.6} in space, we have
\begin{align}\label{GR-KDV_eq:2.7}
 \big({\bm I}+\mathbb{A}^2\big)\frac{d}{dt}{\bm U}+\mathbb{D}(\bm U){\bm U}=0,\ {\bm U}\in\mathbb{V}_h,
\end{align}
with
\begin{align*}
\mathbb{D}({\bm U})=\mathbb{B}+\mathbb{L}_h+\frac{1}{p+2}(\text{diag}\big({\bm U}^p\big)\mathbb{L}_h +\mathbb{L}_h\text{diag}({\bm U}^p)),
\end{align*}
where $\mathbb{A}={\bm I}_{N_{2}}\otimes{\bm D}_{2}^{x}+ {\bm D}_{2}^{y}\otimes {\bm I}_{N_1}$ and
$\mathbb{B}={\bm I}_{N_{2}}\otimes{\bm D}_{3}^{x}+ {\bm D}_{2}^{y}\otimes {\bm D}_{1}^x$. %, and
%$\mathbb{L}_h={\bm I}_{N_{2}}\otimes{\bm D}_{1}^{x}+ {\bm D}_{1}^{y}\otimes {\bm I}_{N_{1}}$.
Note that we have used the equality ${\bm D}_{4}^{w}=({\bm D}_{2}^{w})^2,\ w=x,y$ in the above equation. For more details, please refer to Ref. \cite{GCW14}. In addition, with noting the anti-symmetric property of $\mathbb{B}$ and $\mathbb{L}_h$, we can prove that the matrix ${\mathbb D}(\bm U)$ is anti-symmetric for any mesh function ${\bm U}$.
\begin{lem}\rm\label{GR-KDV-lem2.6}  The semi-discrete system \eqref{GR-KDV_eq:2.7} possesses the following semi-discrete momentum conservation law
\begin{align*}
{P}(t)={P}(0),\ {P}(t)=||{\bm U}||_h^2+|{\bm U}|_{2,h}^2, \ {\bm U}\in\mathbb{V}_h.
\end{align*}

\end{lem}
\begin{prf}\rm Making the discrete inner product of \eqref{GR-KDV_eq:2.7} with ${\bm U}$, we have
\begin{align*}
\frac{d}{dt}\big(||{\bm U}||_h^2+|{\bm U}|_{2,h}^2\big)+\langle \mathbb{D}(\bm U){\bm U},{\bm U}\rangle_h=0.
\end{align*}
With the anti-symmetric property of $\mathbb{D}({\bm U})$, we can obtain
\begin{align*}
\frac{d}{dt}\big(||{\bm U}||_h^2+|{\bm U}|_{2,h}^2\big)=0,
\end{align*}
that is,
\begin{align*}
{P}(t)={P}(0).
\end{align*}
This completes the proof.\qed

\end{prf}

\section{Construction of the linearized Crank-Nicolson momentum-preserving (LCN-MP) scheme} \label{3D_NLS:Sec.2.1}
%In this section, we will propose a LCN-MP scheme by using the linearized Crank-Nicolson method to the semi-discrete system \eqref{GR-KDV_eq:2.7}.
 For a positive integer $M$, let $\Omega_{\tau}=\{t_{n}|t_{n}=n\tau; 0\leq n\leq M\}$
be a uniform partition of $[0,T]$ with time step $\tau=T/M$. Let %$\Omega_{h\tau}=\Omega_{h}\times\Omega_{\tau}$, and denote
$U_{j_1,j_2}^n$ the numerical approximations of $u(x_{j_1},y_{j_2},t_n)$ for $0\le j_r\le N_r-1,r=1,2$ and $0\le n\le M$; denote ${\bm U}^n\in\mathbb{V}_h$ be the solution vector at $t=t_n$ and define
\begin{align*}
&\delta_{t}^{+} {U^{n}_{j_1,j_2}}
=\frac{{ U}_{j_1,j_2}^{n+1}-{U}_{j_1,j_2}^{n}}{\tau},\ { U}_{j_1,j_2}^{n+\frac{1}{2}}=\frac{{U}_{j_1,j_2}^{n+1}+{ U}_{j_1,j_2}^{n}}{2},\
\hat{U}_{j_1,j_2}^{n+\frac{1}{2}}=\frac{3{U}_{j_1,j_2}^{n}-{ U}_{j_1,j_2}^{n-1}}{2},
\end{align*}
for $0\le j_r\le N_r-1,\ r=1,2.$

Applying the linear Crank-Nicolson method to the semi-discrete system \eqref{GR-KDV_eq:2.7} in time, then we can obtain
\begin{align}\label{GR-KDV-LEPSI}
&\big({\bm I}+\mathbb{A}^2\big)\delta_t^+{\bm U}^n+\mathbb{D}(\hat{\bm U}^{n+\frac{1}{2}}){\bm U}^{n+\frac{1}{2}}=0,\ {\bm U}^{n}\in\mathbb{V}_h,\ n=1,\cdots,M-1,
\end{align}
where ${\bm U}^1$ is the solution of the following equation
\begin{align}\label{GR-KDV-LEPSI1}
&\big({\bm I}+\mathbb{A}^2\big)\delta_t^+{\bm U}^0+\mathbb{D}({\bm U}^0){\bm U}^{\frac{1}{2}}=0,\ {\bm U}^{0}\in\mathbb{V}_h,
\end{align}
 which comprises our linearized Crank-Nicolson momentum-preserving (LCN-MP) scheme for the GR-KdV equation. In this paper, for simplicity, we denote $C$ a positive constant which is independent of $h_1$, $h_2$ and $\tau$, and may be different in different case.

%Also we present several discrete Sobolev interpolation formulas and inverse inequalities in the following Lemma. The proof
%can be found in Ref. \cite{zhou90}.
\begin{thm}\rm \label{GR-KDV-thm3.1} The scheme \eqref{GR-KDV-LEPSI}-\eqref{GR-KDV-LEPSI1} possesses the following discrete global momentum conservation law
\begin{align*}
{P}^n=\cdots={P}^0,\ {P}^n=||{\bm U}^n||_h^2+|{\bm U}^n|_{2,h}^2,\ {\bm U}^{n}\in\mathbb{V}_h.
\end{align*}
\end{thm}
\begin{prf}\rm We first show
\begin{align*}
{P}^1={P}^0.
\end{align*}
By nothing
\begin{align*}
\Big\langle\mathbb{D}({\bm U}^0){\bm U}^{\frac{1}{2}},{\bm U}^{\frac{1}{2}}\Big\rangle_h=0,
\end{align*}
we make the discrete inner product of \eqref{GR-KDV-LEPSI1} with ${\bm U}^{\frac{1}{2}}$ and obtain
\begin{align*}
\delta_t^+{P}^0=0,
\end{align*}
that is,
\begin{align*}
{P}^1={P}^0.
\end{align*}
By the similar argument, we have
\begin{align*}
{P}^n=\cdots={P}^1.
\end{align*}
This completes the proof.\qed
\end{prf}

\begin{lem}\rm \label{GR-KDV-lem3.1} Supposing $u_0\in H_p^2(\Omega)$, the solution ${\bm U}^{n}$ of the LCN-MP scheme \eqref{GR-KDV-LEPSI}-\eqref{GR-KDV-LEPSI1} satisfies
\begin{align*}
||{\bm U}^n||_h\le C,\ ||\nabla_h{\bm U}^n||_h\le C,\ ||\Delta_h{\bm U}^n||_{h}\le C,\ ||{\bm U}^n||_{h,\infty}\le C,\ 1\le n\le M.
\end{align*}

\end{lem}
\begin{prf}\rm  According to Theorem \ref{GR-KDV-thm3.1}, we have
 \begin{align}
 ||{\bm U}^n||_h^2+|{\bm U}^n|_{2,h}^2=||{\bm u}_0||_h^2+|{\bm u}_0|_{2,h}^2,\ 1\le n\le M,
 \end{align}
By noting $u_0\in H_p^2(\Omega)$ and Lemma \ref{GR-KDV-lem2.3}, we obtain
 \begin{align}
||{\bm U}^n||_h\le C,\ |{\bm U}^n|_{2,h}\le C,\ 1\le n\le M.
\end{align}
With Lemmas \ref{GR-KDV-lem2.3} and \ref{GR-KDV-lem2.5}, we can get
\begin{align}
||\nabla_h{\bm U}^n||_h\le C,\ ||\Delta_h{\bm U}^n||_{h}\le C,\ ||{\bm U}^n||_{h,\infty}\le C,\ 1\le n\le M.
\end{align}
This completes the proof.\qed
\end{prf}
\begin{thm}\rm The LCN-MP scheme \eqref{GR-KDV-LEPSI}-\eqref{GR-KDV-LEPSI1} is uniquely solvable.

\end{thm}
\begin{prf}\rm For a fixed $n$, the LCN-MP scheme \eqref{GR-KDV-LEPSI}-\eqref{GR-KDV-LEPSI1} can be rewritten as the following linear equation system
\begin{align}
{\bm A}{\bm U}^{n+\frac{1}{2}}={\bm b},\ {\bm U}^{n}\in\mathbb{V}_h,
\end{align}
where ${\bm A}=\big({\bm I}+\mathbb{A}^2+\frac{\tau}{2}\mathbb{D}(\hat{\bm U}^{n+\frac{1}{2}})\big)$ and ${\bm b}=\big({\bm I}+\mathbb{A}^2\big){\bm U}^n.$
In order to obtain the
unique solvability of the scheme, we need to prove that the matrix ${\bm A}$ is invertible.

If ${\bm A}{\bm x}={\bm 0},\ {\bm x}\in\mathbb{V}_h$, we have
\begin{align}
{\bm 0}={\bm x}^T{\bm A}{\bm x}={\bm x}^T\big({\bm I}+\mathbb{A}^2\big){\bm x},
\end{align}
where the anti-symmetry of $\mathbb{D}(\bm U)$ is used. Note that ${\bm I}+\mathbb{A}^2$ is symmetric positive definite, thus,
${\bm x}={\bm 0}$, that is, ${\bm A}{\bm x}={\bm 0}$ has only zero solution. Therefore,  ${\bm A}$ is invertible. This completes  the proof.
\qed
\end{prf}
\begin{lem}\rm\label{GR-KDV-lem3.2}
Let
\begin{align}\label{GR-KDVeq:3.8}
&||{\bm U}^n||_h\le C,\ ||\nabla_h{\bm U}^n||_h\le C,\ ||\Delta_h{\bm U}^n||_{h}\le C,\ ||{\bm U}^n||_{h,\infty}\le C,\ ||{\bm V}^n||_h\le C,\\\label{GR-KDVeq:3.9}
&||\nabla_h{\bm V}^n||_h\le C,\ ||\Delta_h{\bm V}^n||_{h}\le C,\ ||{\bm V}^n||_{h,\infty}\le C,\ {\bm U}^n,{\bm V}^n\in\mathbb{V}_h,
\end{align} for $1\le n\le M$, we then have
\begin{align*}
\Big\langle\mathbb{D}({\bm U}^{0}){\bm U}^{\frac{1}{2}}&-\mathbb{D}({\bm V}^{0}){\bm V}^{\frac{1}{2}},{\bm \eta}^{\frac{1}{2}}\Big\rangle_h
\le C\big(||{{\bm \eta}}^{0}||_h^2+||\Delta_h{\bm \eta}^{0}||_{h}^2+||{\bm \eta}^{1}||_{h}^2+||\Delta_h{\bm \eta}^{1}||_{h}^2\big),
\end{align*}
and
\begin{align*}
\Big\langle\mathbb{D}(\hat{\bm U}^{n+\frac{1}{2}}){\bm U}^{n+\frac{1}{2}}&-\mathbb{D}(\hat{\bm V}^{n+\frac{1}{2}}){\bm V}^{n+\frac{1}{2}},{\bm \eta}^{n+\frac{1}{2}}\Big\rangle_h\nonumber\\
&\le C\big(||{{\bm \eta}}^{n-1}||_h^2+||{{\bm \eta}}^{n}||_h^2+||{{\bm \eta}}^{n+1}||_h^2+||\Delta_h{\bm \eta}^{n}||_{h}^2+||\Delta_h{\bm \eta}^{n+1}||_{h}^2\big),
\end{align*}
where ${\bm \eta}^n={\bm U}^n-{\bm V}^n$ and $n=1,\cdots,M-1$.
\end{lem}
\begin{prf}\rm
Denoting
\begin{align*}
F({\bm x}_1,{\bm x}_2,{\bm x}_3)
&=\frac{1}{p+2}\Big(\frac{3{\bm x}_2-{\bm x}_1}{2}\Big)^p\cdot\mathbb{L}_h\frac{{\bm x}_3+{\bm x}_2}{2}+\frac{1}{p+2}\mathbb{L}_h\Big(\big(\frac{3{\bm x}_2-{\bm x}_1}{2}\big)^p\cdot\frac{{\bm x}_3+{\bm x}_2}{2}\Big)\\
&:=F_1({\bm x}_1,{\bm x}_2,{\bm x}_3)+ F_2({\bm x}_1,{\bm x}_2,{\bm x}_3),\ {\bm x}_i\in\mathbb{V}_h,\ i=1,2,3,
\end{align*}
we have
\begin{small}
\begin{align}\label{GR-KDVeq:3.10}
{\bm f}&:=F_1({\bm U}^{n-1},{\bm U}^{n},{\bm U}^{n+1})-F_1({\bm V}^{n-1},{\bm V}^{n},{\bm V}^{n+1})\nonumber\\
&=\frac{1}{p+2}\Big(\hat{{\bm U}}^{n+\frac{1}{2}}\Big)^p\cdot\mathbb{L}_h{\bm U}^{n+\frac{1}{2}}-\frac{1}{p+2}\Big(\hat{{\bm V}}^{n+\frac{1}{2}}\Big)^p\cdot\mathbb{L}_h{\bm V}^{n+\frac{1}{2}}\nonumber\\
&=\frac{1}{p+2}\Big[\Big(\hat{{\bm U}}^{n+\frac{1}{2}}\Big)^p-\Big(\hat{{\bm V}}^{n+\frac{1}{2}}\Big)^p\Big]\cdot\mathbb{L}_h{\bm U}^{n+\frac{1}{2}}+\frac{1}{p+2}\Big(\hat{{\bm V}}^{n+\frac{1}{2}}\Big)^p\cdot\mathbb{L}_h{\bm \eta}^{n+\frac{1}{2}}\nonumber\\
&=\frac{1}{p+2}\Big[\hat{{\bm \eta}}^{n+\frac{1}{2}}\cdot\sum_{l=0}^{p-1}\Big((\hat{{\bm U}}^{n+\frac{1}{2}})^{p-l-1}\cdot(\hat{{\bm V}}^{n+\frac{1}{2}})^{l}\Big)\Big]\cdot\mathbb{L}_h{\bm U}^{n+\frac{1}{2}}\nonumber\\
&~~~~~~~~~~~~~~~~~~~~~~~~~~~~~~~~~~~~~~~~~~~~~~~~~~~~~~~~~~~~~~+\frac{1}{p+2}\Big(\hat{{\bm V}}^{n+\frac{1}{2}}\Big)^p\cdot\mathbb{L}_h{\bm \eta}^{n+\frac{1}{2}},
\end{align}
\end{small}
and
\begin{small}
\begin{align}\label{GR-KDVeq:3.11}
{\bm g}&:=F_2({\bm U}^{n-1},{\bm U}^{n},{\bm U}^{n+1})-F_2({\bm V}^{n-1},{\bm V}^{n},{\bm V}^{n+1})\nonumber\\
&=\frac{1}{p+2}\mathbb{L}_h\Big(\big(\hat{\bm U}^{n+\frac{1}{2}}\big)^p\cdot{\bm U}^{n+\frac{1}{2}}\Big)-\frac{1}{p+2}\mathbb{L}_h\Big(\big(\hat{\bm V}^{n+\frac{1}{2}}\big)^p\cdot{\bm V}^{n+\frac{1}{2}}\Big)\nonumber\\
&=\frac{1}{p+2}\mathbb{L}_h\Big[\Big(\big(\hat{\bm U}^{n+\frac{1}{2}}\big)^p-\big(\hat{\bm V}^{n+\frac{1}{2}}\big)^p\Big)\cdot{\bm U}^{n+\frac{1}{2}}\Big]+\frac{1}{p+2}\mathbb{L}_h\Big(\big(\hat{\bm V}^{n+\frac{1}{2}}\big)^p\cdot{\bm \eta}^{n+\frac{1}{2}}\Big)\nonumber\\
&=\frac{1}{p+2}\mathbb{L}_h\Big[\hat{\bm \eta}^{n+\frac{1}{2}}\cdot\sum_{l=0}^{p-1}\Big(\big(\hat{\bm U}^{n+\frac{1}{2}}\big)^{p-l-1}\cdot\big(\hat{\bm V}^{n+\frac{1}{2}}\big)^{l}\Big)\cdot{\bm U}^{n+\frac{1}{2}}\Big]\nonumber\\
&~~~~~~~~~~~~~~~~~~~~~~~~~~~~~~~~~~~~~~~~~~~~~~~~~~~~~~~~~~+\frac{1}{p+2}\mathbb{L}_h\Big(\big(\hat{\bm V}^{n+\frac{1}{2}}\big)^p\cdot{\bm \eta}^{n+\frac{1}{2}}\Big).
\end{align}
\end{small}
With noting
\begin{align*}
\Big\langle\big(\mathbb{B}+\mathbb{L}_h\big){\bm \eta}^{n+\frac{1}{2}},{\bm \eta}^{n+\frac{1}{2}}\Big\rangle_h=0,
\end{align*}
we then obtain
\begin{align}\label{GR-KDV:eq:3.8}
\Big\langle\mathbb{D}(\hat{\bm U}^{n+\frac{1}{2}}){\bm U}^{n+\frac{1}{2}}&-\mathbb{D}(\hat{\bm V}^{n+\frac{1}{2}}){\bm V}^{n+\frac{1}{2}},{\bm \eta}^{n+\frac{1}{2}}\Big\rangle_h=\Big\langle{\bm f},{\bm \eta}^{n+\frac{1}{2}}\Big\rangle_h+\Big\langle{\bm g},{\bm \eta}^{n+\frac{1}{2}}\Big\rangle_h.
\end{align}
%and
%\begin{align}\label{GR-KDV:eq:3.9}
%\langle\mathbb{D}(\hat{\bm U}^{n+\frac{1}{2}}){\bm U}^{n+\frac{1}{2}}&-\mathbb{D}(\hat{\bm V}^{n+\frac{1}{2}}){\bm V}^{n+\frac{1}{2}},{\bm \eta}^{n+\frac{1}{2}}\rangle_h=\langle{\bm g},{\bm \eta}^{n+\frac{1}{2}}\rangle_h.
%\end{align}
By using Lemmas \ref{GR-KDV-lem2.1} and \ref{GR-KDV-lem2.3}-\ref{GR-KDV-lem2.5}, and Eqs. \eqref{GR-KDVeq:3.8}-\eqref{GR-KDVeq:3.9}, we can deduce from \eqref{GR-KDVeq:3.10}-\eqref{GR-KDVeq:3.11} that
\begin{small}
\begin{align}\label{GR-KDV:eq:3.10}
\Big\langle{\bm f},{\bm \eta}^{n+\frac{1}{2}}\Big\rangle_h&=\frac{h_1h_2}{p+2}\sum_{j_1=0}^{N_1-1}\sum_{j_2=0}^{N_2-1}\Bigg\{\Big[\hat{{\eta}}_{j_1,j_2}^{n+\frac{1}{2}}\sum_{l=0}^{p-1}\Big((\hat{{ U}}_{j_1,j_2}^{n+\frac{1}{2}})^{p-l-1}\cdot(\hat{{ V}}_{j_1,j_2}^{n+\frac{1}{2}})^{l}\Big)\Big](\mathbb{L}_h{\bm U}^{n+\frac{1}{2}})_{j_1,j_2}\nonumber\\
&-\Big(\hat{{ V}}_{j_1,j_2}^{n+\frac{1}{2}}\Big)^p(\mathbb{L}_h{ \bm \eta}^{n+\frac{1}{2}})_{j_1,j_2}\Bigg\}{ \eta}_{j_1,j_2}^{n+\frac{1}{2}}\nonumber\\
&\le C\big(||\hat{{\bm \eta}}^{n+\frac{1}{2}}||_h^2+||{{\bm \eta}}^{n+\frac{1}{2}}||_{h,\infty}^2||\mathbb{L}_h{\bm U}^{n+\frac{1}{2}}||_{h}^2+
||\mathbb{L}_h{\bm \eta}^{n+\frac{1}{2}}||_{h}^2+||{\bm \eta}^{n+\frac{1}{2}}||_{h}^2\big)\nonumber\\
&\le C\big(||\hat{{\bm \eta}}^{n+\frac{1}{2}}||_h^2+||\Delta_h{\bm \eta}^{n+\frac{1}{2}}||_{h}^2+
||\nabla_h{\bm \eta}^{n+\frac{1}{2}}||_{h}^2+||{\bm \eta}^{n+\frac{1}{2}}||_{h}^2\big)\nonumber\\
&\le C\big(||\hat{{\bm \eta}}^{n+\frac{1}{2}}||_h^2+||\Delta_h{\bm \eta}^{n+\frac{1}{2}}||_{h}^2+||{\bm \eta}^{n+\frac{1}{2}}||_{h}^2\big),\nonumber\\
& \le C\big(||{{\bm \eta}}^{n-1}||_h^2+||{{\bm \eta}}^{n}||_h^2+||{{\bm \eta}}^{n+1}||_h^2+||\Delta_h{\bm \eta}^{n}||_{h}^2+||\Delta_h{\bm \eta}^{n+1}||_{h}^2\big),
\end{align}
\end{small}
and
\begin{small}
\begin{align}\label{GR-KDV:eq:3.11}
\Big\langle{\bm g},{\bm \eta}^{n+\frac{1}{2}}\Big\rangle_h&=\frac{1}{p+2}\Big\langle\hat{\bm \eta}^{n+\frac{1}{2}}\cdot\sum_{l=0}^{p-1}\Big(\big(\hat{\bm U}^{n+\frac{1}{2}}\big)^{p-l-1}\cdot\big(\hat{\bm V}^{n+\frac{1}{2}}\big)^{l}\Big)\cdot{\bm U}^{n+\frac{1}{2}},-\mathbb{L}_h{\bm \eta}^{n+\frac{1}{2}}\Big\rangle_h\nonumber\\
&+\Big\langle\big(\hat{\bm V}^{n+\frac{1}{2}}\big)^p\cdot{\bm \eta}^{n+\frac{1}{2}},-\mathbb{L}_h{\bm \eta}^{n+\frac{1}{2}}\Big\rangle_h\nonumber\\
&\le C\big(||\nabla_h{\bm \eta}^{n+\frac{1}{2}}||_{h}^2+||{\bm \eta}^{n+\frac{1}{2}}||_{h}^2+||\hat{\bm \eta}^{n+\frac{1}{2}}||_{h}^2\big)\nonumber\\
&\le C\big(||\Delta_h{\bm \eta}^{n+\frac{1}{2}}||_{h}^2+||{\bm \eta}^{n+\frac{1}{2}}||_{h}^2+||\hat{\bm \eta}^{n+\frac{1}{2}}||_{h}^2\big)\nonumber\\
&\le  C\big(||{{\bm \eta}}^{n-1}||_h^2+||{{\bm \eta}}^{n}||_h^2+||{{\bm \eta}}^{n+1}||_h^2+||\Delta_h{\bm \eta}^{n}||_{h}^2+||\Delta_h{\bm \eta}^{n+1}||_{h}^2\big).
\end{align}
\end{small}
Substituting \eqref{GR-KDV:eq:3.10} and \eqref{GR-KDV:eq:3.11} into \eqref{GR-KDV:eq:3.8},
we have
\begin{align*}
\Big\langle\mathbb{D}(\hat{\bm U}^{n+\frac{1}{2}}){\bm U}^{n+\frac{1}{2}}&-\mathbb{D}(\hat{\bm V}^{n+\frac{1}{2}}){\bm V}^{n+\frac{1}{2}},{\bm \eta}^{n+\frac{1}{2}}\Big\rangle_h\nonumber\\
&\le C\big(||{{\bm \eta}}^{n-1}||_h^2+||{{\bm \eta}}^{n}||_h^2+||{{\bm \eta}}^{n+1}||_h^2+||\Delta_h{\bm \eta}^{n}||_{h}^2+||\Delta_h{\bm \eta}^{n+1}||_{h}^2\big).
\end{align*}
Similarly, we have
\begin{align*}
\Big\langle\mathbb{D}({\bm U}^{0}){\bm U}^{\frac{1}{2}}&-\mathbb{D}({\bm V}^{0}){\bm V}^{\frac{1}{2}},{\bm \eta}^{\frac{1}{2}}\Big\rangle_h
\le C\big(||{{\bm \eta}}^{0}||_h^2+||\Delta_h{\bm \eta}^{0}||_{h}^2+||{\bm \eta}^{1}||_{h}^2+||\Delta_h{\bm \eta}^{1}||_{h}^2\big).
\end{align*}
\qed

\end{prf}

\section{An a priori estimate}\label{Sec.4}
In this section, we will establish an a priori estimate for the proposed scheme \eqref{GR-KDV-LEPSI}-\eqref{GR-KDV-LEPSI1} in discrete $L^{\infty}$-norm.
For simplicity, we let $\Omega=[0,2\pi]^{2}$. More general cuboid domain can be translated into $\Omega$.
We assume that $C_{p}^{\infty}(\Omega)$ be a set of infinitely differentiable functions with the period $2\pi$ defined on $\Omega$ for all
variables. $H_{p}^{s}(\Omega)$ is the closure of $C_{p}^{\infty}(\Omega)$ in $H^{s}(\Omega)$.
The semi-norm and the norm of $H_{p}^{s}(\Omega)$ are denoted by $\arrowvert\cdot\arrowvert_{s}$
and $\Arrowvert\cdot\Arrowvert_{s}$ respectively.
$\Arrowvert\cdot\Arrowvert_{0}$ is denoted by $\Arrowvert\cdot\Arrowvert$ for simplicity.

Let $N_1=N_2=N$, the interpolation space $S_{N}^{''}$ can be rewritten as
\begin{align*}
&S_{N}^{''}=\Big\{ u| u=\sum_{\arrowvert j_1\arrowvert,\arrowvert j_2\arrowvert\leq\frac{N}{2}}
\frac{{\hat u}_{j_1,j_2}}{c_{j_1}c_{j_2}}e^{\text{i}(j_1x+j_2y)}: {\hat u}_{\frac{N}{2},j_2}={\hat u}_{-\frac{N}{2},j_2},\ {\hat u}_{j_1,\frac{N}{2}}={\hat u}_{j_1,-\frac{N}{2}}\Big\},
\end{align*}
where $c_{l}=1,\ |l|<\frac{N}{2},\ c_{-\frac{N}{2}}=c_{\frac{N}{2}}=2$. The projection space is defined as
\begin{align*}
S_{N}=\Big\{u| u=\sum_{\arrowvert j_1\arrowvert,\arrowvert j_2\arrowvert\leq\frac{N}{2}}
{\tilde u}_{j_1,j_2}e^{\text{i}(j_1x+j_2y)}\Big\}.
\end{align*}
It is clear to see that $S_{N-2}\subseteq S_{N}^{''}\subseteq S_{N}$.
We denote by $P_{N}:L^{2}(\Omega)\to S_{N}$ as the orthogonal projection operator
and recall the interpolation operator $I_{N}:C(\Omega)\to S_{N}^{''}$.
Further, $P_{N}$ and $I_{N}$ satisfy \cite{GWWC17}:
\begin{align*}
&1.\  P_{N}\partial_{w}u=\partial_{w} P_{N}u,\  I_{N}\partial_{w}u\ne \partial_{w} I_{N}u,\ w=x,\ \text{or}\ y. \\
&2.\  P_{N}u=u,\ \forall u\in S_{N}, \ I_{N}u=u,\ \forall u\in S_{N}^{''}.
\end{align*}
%Now, we introduce some useful Lemmas that will be used frequently in the subsequent error estimate.
\begin{lem}\rm \label{GR-KDV:lem4.1} \cite{GWWC17}
 For ${u}\in S_{N}^{''}$,
$\Arrowvert { u}\Arrowvert\leq \Arrowvert {\bm u}\Arrowvert_{h}\leq 2\Arrowvert { u}\Arrowvert$, where ${\bm u}\in\mathbb{V}_h$.
\end{lem}
\begin{lem}\rm \label{GR-KDV:lem4.2}
 \cite{CQ82} If $0\leq l \leq s$ and ${u}\in H_{p}^{s}(\Omega)$, then
\begin{align*}
&\Arrowvert P_{N}{u}-{u}\Arrowvert_{l}\leq CN^{l-s}\arrowvert {u}\arrowvert_{s},\\
&||P_{N} u||_{l}\le C ||u||_{l},
\end{align*}
 and in addition if $s>1$ then
\begin{align*}
&\Arrowvert I_{N}{ u}-{ u}\Arrowvert_{l}\leq CN^{l-s}\arrowvert {u}\arrowvert_{s}.
%&||I_{N} u||_{l}\le C ||u||_{l}.
\end{align*}
\end{lem}
\begin{lem}\rm \label{GR-KDV:lem4.3} \cite{GWWC17} For ${ u}\in H_{p}^{s}(\Omega),\ s>1$, let ${u}^{*}=P_{N-2}{u},\ N>2$. Then, we have
 \begin{align*}%\label{3d_DNLS:4.5}
 \Arrowvert {\bm u}^{*}-{\bm u}\Arrowvert_{h}\leq CN^{-s} \arrowvert {u}\arrowvert_{s},\ {\bm u}^{*},{\bm u}\in\mathbb{V}_h.
 \end{align*}
\end{lem}
\begin{lem}\rm \cite{GW12}\label{GR-KDV:lem4.4} For any ${ u}\in S_{qN}$, we have
 \begin{align*}%\label{3d_DNLS:4.5}
 \Arrowvert I_N u\Arrowvert_{l}\leq q\Arrowvert u \Arrowvert_{l}.
 \end{align*}
\end{lem}
\begin{lem}\rm \label{GR-KDV:lem4.5} For ${ u}\in H_{p}^{s+1}(\Omega),\ s>1$, let ${u}^{*}=P_{N-2}{u},\ N>2$. Then, we have
 \begin{align*}%\label{3d_DNLS:4.5}
 \Arrowvert\nabla_h\big( {\bm u}^{*}-{\bm u}\big)\Arrowvert_{h}\leq CN^{-s} \arrowvert {u}\arrowvert_{s+1},\ {\bm u}^{*},{\bm u}\in\mathbb{V}_h.
 \end{align*}
\end{lem}
The proof is similar to Lemma 4.4 in Ref. \cite{JCW17}. For brevity, we omit it.
\begin{lem}\rm \label{GR-KDV:lem4.6} For ${ u}\in H_{p}^{s+2}(\Omega),\ s>1$, let ${u}^{*}=P_{N-2}{u},\ N>2$. Then, we have
 \begin{align*}%\label{3d_DNLS:4.5}
 \Arrowvert\Delta_h\big( {\bm u}^{*}-{\bm u}\big)\Arrowvert_{h}\leq CN^{-s} \arrowvert {u}\arrowvert_{s+2},\ {\bm u}^{*},{\bm u}\in\mathbb{V}_h.
 \end{align*}
\end{lem}
\begin{prf}\rm According to Lemmas \ref{GR-KDV-lem2.3} and \ref{GR-KDV:lem4.1}, we have
\begin{align}\label{GR-KDV:eq:4.1}
\Arrowvert\Delta_h\big( {\bm u}^{*}-{\bm u}\big)\Arrowvert_{h}\le |{\bm u}^{*}-{\bm u}|_{2,h}&=\Arrowvert\mathbb{A}({\bm u}^{*}-{\bm u})\Arrowvert_h\le 2\Arrowvert I_N(\Delta(I_N({u}^{*}-{u})))\Arrowvert,
\end{align}
where we have used the fact $\big[\mathbb{A}({\bm u}^{*}-{\bm u})\big]_{j_1,j_2}=\big[I_N(\Delta(I_N({u}^{*}-{u})))\big](x_{j_1},y_{j_2})$.
By noting $\Delta(I_N({u}^{*}-{u}))\in S_{2N}$, we can deduce from Lemmas \ref{GR-KDV:lem4.2} and \ref{GR-KDV:lem4.4} that
\begin{align}\label{GR-KDV:eq:4.2}
\Arrowvert I_N(\Delta(I_N({u}^{*}-{u})))\Arrowvert&\leq 2\Arrowvert\Delta(I_N({u}^{*}-{ u}))\Arrowvert\nonumber\\
&\le 2\big(\Arrowvert\Delta(u^{*}-u)\Arrowvert+\Arrowvert\Delta (u-I_Nu)\Arrowvert\big)\nonumber\\
&\le C\big(\Arrowvert u^{*}-u\Arrowvert_2+\Arrowvert u-I_Nu\Arrowvert_2\big)\nonumber\\
&\le CN^{-s}|u|_{s+2}.
\end{align}
Substituting \eqref{GR-KDV:eq:4.2} into \eqref{GR-KDV:eq:4.1}, we finish the proof.\qed

\end{prf}

\begin{lem}\label{GR-KDV:lem4.9}\rm Let $u^*=P_{N-2}u$ and
\begin{align*}
 u_t^*+\Delta^2u_t^*&+\Delta u_x^*+\mathbb{L} u^*+\frac{1}{p+2}((u^*)^p\mathbb{L}u^* +\mathbb{L}(I_N(u^*)^{p+1}))=\nonumber\\
 &u_t+\Delta^2u_t+\Delta u_x+\mathbb{L} u+\frac{1}{p+2}(u^p\mathbb{L}u +\mathbb{L}(u^{p+1}))+\xi_1.
\end{align*}
If $u\in C^1(0,T; H^{s+4}(\Omega)),\ s>1$, we have
\begin{align*}
||{\bm \xi}_1||_h\le CN^{-s},\ {\bm \xi}_1\in \mathbb{V}_{h}.
 \end{align*}

\end{lem}
\begin{prf}\rm
With Lemma \ref{GR-KDV:lem4.2}, we can obtain the following approximation estimate
\begin{align}\label{GR-KDV:eq:4.4}
||\partial_t^k(u^*-u)||_{l}=||P_{N-2}(\partial_t^ku)-\partial_t^ku||_{l}\le CN^{-s}|\partial_t^ku|_{s+l}.
\end{align}
Let
\begin{align}\label{GR-KDV:eq:4.5}
\tilde{\xi}&=((u^*)^p\mathbb{L}u^* +\mathbb{L}(I_N(u^*)^{p+1}))-(u^p\mathbb{L}u +\mathbb{L}(u^{p+1}))\nonumber\\
&=\Big[((u^*)^p-u^p)\mathbb{L}u^*+u^p\mathbb{L}(u^*-u)\Big]\nonumber\\
&~~~~+\Big[\mathbb{L}(I_N(u^*)^{p+1}-(u^*)^{p+1})+\mathbb{L}((u^*)^{p+1}-u^{p+1})\Big]\nonumber\\
&=:\tilde{\xi}_1+\tilde{\xi}_2.
\end{align}
With Lemma \ref{GR-KDV:lem4.2}, we have
\begin{align}\label{GR-KDV:eq:4.6}
||\tilde{\xi}_1||&\le ||((u^*)^p-u^p)\mathbb{L}u^*||+||u^p\mathbb{L}(u^*-u)||\nonumber\\
%&\le ||\mathbb{L}u^*||_{\infty}||((u^*)^p-u^p)||+||u^p||_{\infty}||\mathbb{L}(u^*-u)||\nonumber\\
&\le ||\mathbb{L}u^*||_{L^{\infty}}||(u^*)^p-u^p||+||u||_{L^{\infty}}^p||\mathbb{L}(u^*-u)||\nonumber\\
& \le ||\mathbb{L}u^*||_{L^{\infty}}\Big|\Big|\sum_{l=0}^{p-1}(u^*)^{p-l-1}u^l\Big|\Big|_{L^{\infty}}||u^*-u||+||u||_{L^{\infty}}^p||\mathbb{L}(u^*-u)||\nonumber\\
&\le C||u^*-u||_1\nonumber\\
&\le CN^{-s},
\end{align}
and
\begin{align}\label{GR-KDV:eq:4.7}
||\tilde{\xi}_2||&\le ||\mathbb{L}(I_N(u^*)^{p+1}-(u^*)^{p+1}))||+||\mathbb{L}((u^*)^{p+1}-u^{p+1})||\nonumber\\
&\le C\Big( ||I_N(u^*)^{p+1}-(u^*)^{p+1}||_1+||(u^*)^{p}(\mathbb{L}u^{*}-\mathbb{L}u)||+||((u^*)^{p}-u^{p})\mathbb{L}u||\Big)\nonumber\\
&\le C\Big( ||I_N(u^*)^{p+1}-(u^*)^{p+1}||_1+||(u^*)^{p}(\mathbb{L}u^{*}-\mathbb{L}u)||+||((u^*)^{p}-u^{p})\mathbb{L}u||\Big)\nonumber\\
&\le  CN^{-s}|(u^*)^{p+1}|_{s+1}+ ||u^*||_{L^{\infty}}^p||\mathbb{L}u^{*}-\mathbb{L}u||\nonumber\\
&~~~~~~~~+||\mathbb{L}u||_{L^{\infty}}||\sum_{l=0}^{p-1}((u^*)^{p-l-1})u^l||_{L^{\infty}}||u^*-u||\nonumber\\
&%\le  CN^{-s}|(u^*)^{p+1}|_{s+1}+ CN^{-s}|u|_{s+1}.
\le  CN^{-s}.
\end{align}
Thus, we can deduce from \eqref{GR-KDV:eq:4.6} and \eqref{GR-KDV:eq:4.7} that
\begin{align}\label{GR-KDV:eq:4.8}
||\tilde{\xi}||\le ||\tilde{\xi}_1||+||\tilde{\xi}_2||\le CN^{-s}.%|(u^*)^{p+1}|_{s+1}+ CN^{-s}|u|_{s+1}.
\end{align}
With \eqref{GR-KDV:eq:4.4} and \eqref{GR-KDV:eq:4.8}, we have
\begin{align}\label{GR-KDV:eq:4.9}
||\xi_1||\le CN^{-s}.
\end{align}
Further, by noting $\xi_1\in S_{(p+1)N}$, we can obtain
\begin{align}\label{GR-KDV:eq:4.9n}
||{\bm \xi}_1||_h\le 2||I_N\xi_1||\le C||\xi_1||\le CN^{-s},
\end{align}
where Lemmas \ref{GR-KDV:lem4.1} and \ref{GR-KDV:lem4.4} are used.\qed
\end{prf}
\begin{lem}\rm \label{GR-KDV-lem4.8}
Let
\begin{align}\label{GR-KDV:eq:4.10}
 ({\bm I}+\mathbb{A}^2)\delta_t^+({\bm u}^*)^0&+\mathbb{D}(({\bm u}^*)^{0})({\bm u}^*)^{\frac{1}{2}}={\bm \xi}^{0},\ ({\bm u}^*)^{0}, {\bm \xi}^0\in\mathbb{V}_h,
\end{align}
and
\begin{align}\label{GR-KDV:eq:4.11}
 &({\bm I}+\mathbb{A}^2)\delta_t^+({\bm u}^*)^n+\mathbb{D}((\hat{\bm u}^*)^{n+\frac{1}{2}})({\bm u}^*)^{n+\frac{1}{2}}={\bm \xi}^{n+\frac{1}{2}},\ ({\bm u}^*)^n,{\bm \xi}^{n+\frac{1}{2}}\in\mathbb{V}_h,
\end{align}
for $n=1,2,\cdots,M-1$.
 If $u\in C^{3}\Big(0,T; H_{p}^{s+4}(\Omega)\Big),\ s>1$, we then have
\begin{align*}\label{GR-KDV:eq:4.12}
||{\bm \xi}^{0}||_h\le C(N^{-s}+\tau),\ ||{\bm \xi}^{n+\frac{1}{2}}||_h\le C(N^{-s}+\tau^2).
\end{align*}
\begin{prf}\rm We denote
\begin{align*}
(1+\Delta^2)\delta_t^+({u}^*)^0&+\Delta({u}^*)_x^{\frac{1}{2}}+\mathbb{L} ({u}^*)^{\frac{1}{2}}\nonumber\\
&=(1+\Delta^2) \partial_{t}{u}^*(x,y,0)+(\Delta\partial_x+\mathbb{L}){u}^*(x,y,0)+{\xi}_2^{0},
\end{align*}
and
\begin{align*}
\frac{1}{p+2}\Big[&(({u}^*)^{0})^p\mathbb{L}({u}^*)^{\frac{1}{2}}+\mathbb{L}I_N\Big((({u}^*)^{0})^{p}({u}^*)^{\frac{1}{2}}\Big)\Big]\nonumber\\
&=\frac{1}{p+2}\Big[({u}^*(x,y,0))^p\mathbb{L}{u}^*(x,y,0)+\mathbb{L}I_N\Big(({u}^*(x,y,0))^{p}{u}^*(x,y,0)\Big)\Big]+{\xi}_3^{0}.
\end{align*}
By the Taylor expansion, we have
\begin{align}
&\delta_t^+({u}^*)^0=\partial_t u^*(x,y,0)+\tau c_1\partial_{tt}u^{*}(x,y,\zeta_1),\ \zeta_1\in(0,\tau),\\
&({u}^*)^{\frac{1}{2}}=u^*(x,y,0)+\tau c_2\partial_{t}u^{*}(x,y,\zeta_2),\ \zeta_2\in(0,\tau),
\end{align}
where $c_1$ and $c_2$ are constants. With noting $u\in C^{3}\Big(0,T; H_{p}^{s+4}(\Omega)\Big),\ s>1$, we have
\begin{align*}
||{\xi}_2^{0}||&\le C\tau\Big(||(1+\Delta^2)\partial_{tt}{u}^*(x,y,\zeta_1)||+||(\Delta\partial_x+\mathbb{L})\partial_t{u}^*(x,y,\zeta_2)||\Big)\nonumber\\
%&\le C\tau(||\Delta^2\partial_{tt}u^*(0)||+||(\Delta\partial_x+\mathbb{L})\partial_tu^*(0)||)\nonumber\\
&\le C\tau (||\partial_{tt}u^*(x,y,\zeta_1)||_4+||\partial_{t}u^*(x,y,\zeta_2)||_3)\nonumber\\
&\le C\tau(||\partial_{tt}u(x,y,\zeta_1)||_4+||\partial_{t}u(x,y,\zeta_2)||_3)\nonumber\\
&\le C\tau,
\end{align*}
and
\begin{align}\label{GR-KDV:4.15}
||{\xi}_3^{0}||&%\le C\tau\big(||({ u}^*(0))^p\mathbb{L}\partial_t{ u}^*(0)||_h+||\mathbb{L}\Big(({u}^*(0))^{p}\partial_t{ u}^*(0)\Big)||_h\big)\nonumber\\
\le C\tau\Big(||({u}^*(x,y,0))^p\cdot\mathbb{L}\partial_t{u}^*(x,y,\zeta_2)||+||\mathbb{L}I_N\Big(({ u}^*(x,y,0))^{p}\cdot\partial_t{u}^*(x,y,\zeta_2)\Big)||\Big)\nonumber\\
&\le C\tau\Big(||({u}^*(x,y,0))^p\cdot\mathbb{L}\partial_t{u}^*(x,y,\zeta_2)||+||\mathbb{L}\Big(({ u}^*(x,y,0))^{p}\cdot\partial_t{u}^*(x,y,\zeta_2)\Big)||\Big)\nonumber\\
&\le C\tau\Big(||{u}^*(x,y,0)||_{L^{\infty}}^p\cdot||\mathbb{L}\partial_t{u}(x,y,\zeta_2)^*||+||\mathbb{L}({ u}^*(x,y,0))^{p}||_{L^{\infty}}\cdot||\partial_t{u}^*(x,y,\zeta_2)||\nonumber\\
&~~~+||\mathbb{L}\partial_t{u}^*(x,y,\zeta_2)||\cdot||{ u}^*(x,y,0)||_{L^{\infty}}^{p}\Big)\nonumber\\
%&\le C\tau\Big(||{u}(0)||_{2}^p\cdot||\partial_t{u}(\zeta_2)||_1+||\partial_t{u}(\zeta_2)||_1\Big)\nonumber\\
&\le C\tau.
\end{align}

Let
\begin{align}\label{GR-KDV-4.16}
(1+\Delta^2)\delta_t^+({u}^*)^n&+\Delta\partial_x({u}^*)^{n+\frac{1}{2}}+\mathbb{L} ({u}^*)^{n+\frac{1}{2}}\nonumber\\
&=(1+\Delta^2)\partial_{t}{ u}^*(x,y,t_{n+\frac{1}{2}})+(\Delta\partial_x+\mathbb{L}){u}^*(x,y,t_{n+\frac{1}{2}})+{\xi}_2^{n+\frac{1}{2}},
\end{align}
and
\begin{align}\label{GR-KDV-4.17}
&\frac{1}{p+2}\Big[((\hat{u}^*)^{n+\frac{1}{2}})^p\mathbb{L}({u}^*)^{n+\frac{1}{2}}+\mathbb{L}I_N\Big(((\hat{u}^*)^{n+\frac{1}{2}})^{p}({ u}^*)^{n+\frac{1}{2}}\Big)\Big]\nonumber\\
&~~~~~~~~~=\frac{1}{p+2}\Big[({u}^*(x,y,t_{n+\frac{1}{2}}))^p\cdot\mathbb{L}{u}^*(x,y,t_{n+\frac{1}{2}})\nonumber\\
&~~~~~~~~~~~~~~~~~~~~~~~~~~+\mathbb{L}I_N\Big(({u}^*(x,y,t_{n+\frac{1}{2}}))^{p}{ u}^*(x,y,t_{n+\frac{1}{2}})\Big)\Big]+{ \xi}_3^{n+\frac{1}{2}}.
\end{align}
By the Taylor expansion, we have
\begin{align}
&\delta_t^+({u}^*)^n=\partial_tu^*(x,y,t_{n+\frac{1}{2}})+c_3\tau^2 \partial_{ttt}u^{*}(x,y,\zeta_3),\ \zeta_3\in(t_n,t_n+\tau),\\
&({u}^*)^{n+\frac{1}{2}}=u^*(x,y,t_{n+\frac{1}{2}})+c_4\tau^2 \partial_{tt}u^{*}(x,y,\zeta_4),\ \zeta_4\in(t_n,t_n+\tau),
\end{align}
where $c_3$ and $c_4$ are constants. An argument similar to \eqref{GR-KDV-4.16} and \eqref{GR-KDV-4.17} used in \eqref{GR-KDV:4.15} shows that
\begin{align}
||{\xi}_2^{n+\frac{1}{2}}||&\le C\tau^2,
\end{align}
and
\begin{align}
||{\xi}_3^{n+\frac{1}{2}}||&\le C\tau^2.%(||\partial_{tt}u(\zeta_4)||_4+||\partial_{tt}u(\zeta_4)||_1)\le C\tau^2, %C\tau\big(||({\bm u}^*(t_{n+\frac{1}{2}}))^p\cdot\mathbb{L}_h\partial_t{\bm u}^*(t_{n+\frac{1}{2}})||_h+||\mathbb{L}_h\Big(({\bm u}^*(t_{n+\frac{1}{2}}))^{p}\cdot\partial_t{\bm u}^*(t_{n+\frac{1}{2}})\Big)||_h\big)\nonumber\\
%&\le C\tau\big(||({u}^*(t_{n+\frac{1}{2}}))^p\cdot\mathbb{L}\partial_t{u}^*(t_{n+\frac{1}{2}})||+||\mathbb{L}I_N\Big(({ u}^*(0))^{p}\cdot\partial_t{u}^*(0)\Big)||\big)\nonumber\\
%&\le C\tau\big(||\partial_t{u}(t_{n+\frac{1}{2}})||_1+||({ u}(t_{n+\frac{1}{2}}))^{p}\cdot\partial_t{u}(t_{n+\frac{1}{2}})||_1\big)\le C\tau,
\end{align}
Noting ${\xi}_2^{0},\ {\xi}_2^{n+\frac{1}{2}}\in S_{2N}$ and ${\xi}_3^{0},\ {\xi}_3^{n+\frac{1}{2}}\in S_{(p+1)N}$, then by using Lemma \ref{GR-KDV:lem4.4}, we can prove
\begin{align}\label{GR-KDV:eq:4.16}
||{\bm \xi}_2^{0}||_h+||{\bm \xi}_3^{0}||_h\le2\big(||I_N{\xi}_2^{0}||+||I_N{ \xi}_3^{0}||\big)\le C\big(||{\xi}_2^{0}||+||{ \xi}_3^{0}||\big)\le C\tau,
\end{align}
and
\begin{align}\label{GR-KDV:eq4.17}
||{\bm \xi}_2^{n+\frac{1}{2}}||_h+||{\bm \xi}_3^{n+\frac{1}{2}}||_h&\le 2\big(||I_N{\xi}_2^{n+\frac{1}{2}}||+||I_N{ \xi}_3^{n+\frac{1}{2}}||\big)\nonumber\\
&\le C\big(||{\xi}_2^{n+\frac{1}{2}}||+||{ \xi}_3^{n+\frac{1}{2}}||\big)\nonumber\\
&\le C\tau^2.
\end{align}
It is clear to see that
\begin{align*}
 &{\bm \xi}^{0}={\bm \xi}_1(0)+{\bm \xi}_2^{0}+{\bm \xi}_3^{0},\ {\bm \xi}^{n+\frac{1}{2}}={\bm \xi}_1(t_{n+\frac{1}{2}})+{\bm \xi}_2^{n+\frac{1}{2}}+{\bm \xi}_3^{n+\frac{1}{2}}.\
\end{align*}
Thus, with Lemma \ref{GR-KDV:lem4.9}, we can deduce from \eqref{GR-KDV:eq:4.16} and \eqref{GR-KDV:eq4.17} that
\begin{align*}
&||{\bm \xi}^{0}||_h\le ||{\bm \xi}_1(0)||_h+||{\bm \xi}_2^{0}||_h+||{\bm \xi}_3^{0}||_h\le C(N^{-s}+\tau),\\
 &||{\bm \xi}^{n+\frac{1}{2}}||_h\le ||{\bm \xi}_1(t_{n+\frac{1}{2}})||_h+||{\bm \xi}_2^{n+\frac{1}{2}}||_h+||{\bm \xi}_3^{n+\frac{1}{2}}||_h\le C(N^{-s}+\tau^2).
 \end{align*}
 This completes the proof.\qed
\end{prf}
\end{lem}

%\begin{lem}\rm \label{GR-KDV:lem5.7}
% (Gronwall inequality \cite{zhou90}).
%Suppose that the discrete function $\big\{\omega^{n}|n=0,1,2,\cdots,M;M\tau=T\big\}$ is nonnegative and satisfies the inequality
%\begin{align*}
%\omega^{n}\leq \widehat{E}+\tau\sum_{l=1}^{n}E_{l}\omega^{l},
%\end{align*}
%where $\widehat{E}$ and $E_{k} (k=1,2,\cdots,M)$ are nonnegative constants. Then
%\begin{align*}
%\max\limits_{0 \leq n\leq
%M}|\omega^{n}|\leq \widehat{E} e^{2\sum_{k=1}^{M}E_{k}\tau},
%\end{align*}
%where $\tau$ is sufficiently small, such that $\tau\big(\max\limits_{k=0,1,\cdots,
%M}E_{k}\big)\leq\frac{1}{2}$.
%\end{lem}

We define the error function by
\begin{align*}
e_{j_1,j_2}^n=(u^*)_{j_1,j_2}^n-U_{j_1,j_2}^n,\ (j_1,j_2)\in J_h^{'},\ 1\le n\le M.
\end{align*}
Subtracting \eqref{GR-KDV-LEPSI} and \eqref{GR-KDV-LEPSI1} from \eqref{GR-KDV:eq:4.11} and \eqref{GR-KDV:eq:4.10}, respectively, we can get
\begin{align}\label{GR-KDV:eq:4.17}
&\big({\bm I}+\mathbb{A}^2\big)\delta_t^+{\bm e}^0+\mathbb{D}(({\bm u}^*)^0)({\bm u}^*)^{\frac{1}{2}}-\mathbb{D}({\bm U}^0){\bm U}^{\frac{1}{2}}={\bm \xi}^0,
\end{align}
and
\begin{align}\label{GR-KDV:eq:4.20}
&\big({\bm I}+\mathbb{A}^2\big)\delta_t^+{\bm e}^n+\mathbb{D}((\hat{\bm u}^*)^{n+\frac{1}{2}})({\bm u}^*)^{n+\frac{1}{2}}-\mathbb{D}(\hat{\bm U}^{n+\frac{1}{2}}){\bm U}^{n+\frac{1}{2}}={\bm \xi}^{n+\frac{1}{2}},
\end{align}
where ${\bm e}^n=({\bm u}^*)^n-{\bm U}^n\in\mathbb{V}_h,\ n=1,2,\cdots,M-1$.

\begin{thm}\rm \label{3d_DNLSlem4.5} We assume $u\in C^{2}\Big(0,T; H_{p}^{s+4}(\Omega)\Big),\ s>1$. Then, there exists a constant $\tau_0>0$ sufficiently small,
 such that, when $0< \tau\leq \tau_0$, we have
\begin{align*}
&||{\bm u}^{1}-{\bm U}^{1}||_{h}+||\Delta_h({\bm u}^{1}-{\bm U}^{1})||_{h}\leq C(N^{-s}+\tau^{2}),
\end{align*}
and
\begin{align*}
||{\bm u}^{1}-{\bm U}^{1}||_{h,\infty}\leq C(N^{-s}+\tau^{2}),\ {\bm u}^{1},{\bm U}^{1}\in\mathbb{V}_h.
\end{align*}
\end{thm}

\begin{prf}\rm  Making the discrete inner product of \eqref{GR-KDV:eq:4.17} with ${\bm e}^{\frac{1}{2}}$, when $0< \tau\leq \tau_0$, we have
\begin{align}\label{GR-KDV:20}
&F^1-F^0\le C\tau\big(F^1+F^0\big)+C(N^{-s}+\tau^{2})^2,
\end{align}
with
\begin{align*}
F^n=||{\bm e}^n||_h^2+|{\bm e}^n|_{2,h}^2,\ n=0,1,
\end{align*}
where Lemmas \ref{GR-KDV-lem3.2} and \ref{GR-KDV-lem4.8} are used.  By Lemmas  \ref{GR-KDV-lem2.3} and \ref{GR-KDV:lem4.3}-\ref{GR-KDV:lem4.6}, we can deduce that
\begin{align}\label{GR-KDV:21}
F^0=||{\bm e}^0||_h^2+|{\bm e}^0|_{2,h}^2&\le ||{\bm e}^0||_h^2+\frac{\pi^4}{16}||\Delta_h{\bm e}^0||_{h}^2\nonumber\\
&=||({\bm u}^*)^0-{\bm u}^0||_h^2+\frac{\pi^4}{16}||\Delta_h(({\bm u}^*)^0-{\bm u}^0)||_h^2\nonumber\\
&\le CN^{-2s}.
\end{align}
With \eqref{GR-KDV:21}, when $\tau$ is sufficiently small, such that, when $0< \tau\leq \tau_0$,  we can get from \eqref{GR-KDV:20} that
\begin{align}\label{GR-KDV-26}
||{\bm e}^1||_h^2+|{\bm e}^1|_{2,h}^2\le C(N^{-s}+\tau^{2})^2.
\end{align}
%which implies that
%\begin{align*}%\label{GR-KDV:eq4.18}
%||{\bm e}^1||_h+|{\bm e}^1|_{2,h}\le C(N^{-s}+\tau^{2}).
%\end{align*}
According to Lemma \ref{GR-KDV-lem2.3}, we get
\begin{align}
||{\bm e}^1||_h^2+||\Delta_h{\bm e}^1||_{h}^2\le C(N^{-s}+\tau^{2})^2,
\end{align}
which implies that
\begin{align}\label{GR-KDV:eq4.18}
||{\bm e}^1||_h+||\Delta_h{\bm e}^1||_{h}\le C(N^{-s}+\tau^{2}).
\end{align}
By using Lemmas \ref{GR-KDV:lem4.3} and \ref{GR-KDV:lem4.6}, and Eq. \eqref{GR-KDV:eq4.18}, we have
\begin{align}\label{GR-KDV-eq4.19}
||{\bm u}^{1}-{\bm U}^{1}||_{h}&+||\Delta_h({\bm u}^{1}-{\bm U}^{1})||_{h}\nonumber\\
&\le \big(||{\bm u}^{1}-({\bm u}^*)^1||_h+||\Delta_h{\bm u}^{1}-\Delta_h({\bm u}^*)^1||_h+||{\bm e}^1||_{h}+||\Delta_h{\bm e}^{1}||_{h}\big)\nonumber\\
&\le C(N^{-s}+\tau^{2}).
\end{align}
With Lemma \ref{GR-KDV-lem2.5}, we can deduce from \eqref{GR-KDV-eq4.19} that
\begin{align*}
||{\bm u}^{1}-{\bm U}^{1}||_{h,\infty}\le C(N^{-s}+\tau^{2}).
\end{align*}
This completes the proof.\qed

%\begin{align}
% \delta_t^+({\bm U}^*)^0+\delta_t^+\mathbb{A}({\bm U}^*)^0&+\mathbb{B}({\bm U}^*)^{\frac{1}{2}}+\mathbb{L}_h ({\bm U}^*)^{\frac{1}{2}}
% +\frac{1}{p+2}(diag(({\bm U}^*)^{0})^p)\mathbb{L}_h({\bm U}^*)^{\frac{1}{2}}\nonumber\\
% &+\mathbb{L}_h diag((({\bm U}^*)^{0})^{p})({\bm U}^*)^{\frac{1}{2}}={\bm \xi}^{0},
%\end{align}
\end{prf}

\begin{thm}\rm \label{3d_DNLSlem4.5} We assume $u\in C^{3}\Big(0,T; H_{p}^{s+4}(\Omega)\Big),\ s>1$. %Then, there exists a constant $\tau_0>0$ sufficiently small,
 %such that, when $0< \tau\leq \tau_0$,
Then, when $\tau$ is sufficiently
small, such that, $C\tau\le \frac{1}{2}$, we have
\begin{align*}
&||{\bm u}^{n}-{\bm U}^{n}||_{h}+||\Delta_h({\bm u}^{n}-{\bm U}^{n})||_{h}\leq C(N^{-s}+\tau^{2}),\  ||{\bm U}^n||_{h,\infty}\le C(N^{-s}+\tau^{2}),
\end{align*}
where ${\bm u}^{n},{\bm U}^{n}\in\mathbb{V}_h,\ n=2,3,\cdots,M$.
\end{thm}
\begin{prf}\rm %\begin{align}
% \delta_t^+(u^*)^n&+\Delta^2\delta_t^+(u^*)^n+\Delta (u_x^*)^{n+\frac{1}{2}}+\mathbb{L} (u^*)^{n+\frac{1}{2}}\nonumber\\
% &+\frac{1}{p+2}(((\hat{u}^*)^{n+\frac{1}{2}})^p\mathbb{L}(u^*)^{n+\frac{1}{2}}
% +\mathbb{L}(I_N((\hat{u}^*)^{n+\frac{1}{2}})^{p}(u^*)^{n+\frac{1}{2}}))=\xi^{n+\frac{1}{2}},
%\end{align}
%Noting $\xi^{n+\frac{1}{2}}\in S_{(p+1)N}$
%\begin{align}
%||{\bm \xi}^{n+\frac{1}{2}}||_h=||I_N\xi^{n+\frac{1}{2}}||_h\le 2||I_N\xi^{n+\frac{1}{2}}||\le 2(p+1)||\xi^{n+\frac{1}{2}}||\le C(N^{-s}+\tau).
%\end{align}
%\begin{align}
%&\big({\bm I}+\mathbb{A}^2\big)\delta_t^+({\bm u}^*)^n+\mathbb{D}(({\bm u}^*)^{n+\frac{1}{2}})({\bm u}^*)^{n\frac{1}{2}}={\bm \xi}^{n+\frac{1}{2}}.
%\end{align}

Making the discrete inner product of \eqref{GR-KDV:eq:4.20} with ${\bm e}^{n+\frac{1}{2}}$, we then have
\begin{align}\label{GR-KDV:eq:4.21}
&F^{n+1}-F^n\le C\tau\big(F^{n+1}+F^n\big)+C\tau||{\bm e}^{n-1}||_h^2+C\tau(N^{-s}+\tau^2)^2,
\end{align}
with
\begin{align*}
F^n=||{\bm e}^n||_h^2+|{\bm e}^n|_{2,h}^2,\ 1\le n\le M,%||{\bm e}^{n-1}||_h^2+||\Delta_h{\bm e}^{n-1}||_h^2+
\end{align*}
where Lemmas \ref{GR-KDV-lem3.2} and \ref{GR-KDV-lem4.8} are used. Summing up for $n$ from 1 to $m$ and then replacing $m$ by $n-1$, we can get from \eqref{GR-KDV:eq:4.21} that
\begin{align}\label{GR-KDV:eq:4.22}
F^{n}&\le F^1+C\tau\sum_{l=1}^nF^{l}+C\tau||{\bm e}^{0}||_h^2+CT(N^{-s}+\tau^2)^2\nonumber\\
&\le C\tau\sum_{l=1}^nF^{l}+CT(N^{-s}+\tau^2)^2,
\end{align}
where \eqref{GR-KDV-26} and \eqref{GR-KDV:eq4.18} are used.
Applying the Gronwall inequality \cite{zhou90} to \eqref{GR-KDV:eq:4.22}, then we have
\begin{align}
||{\bm e}^n||_h^2+|{\bm e}^n|_{2,h}^2\le C(N^{-s}+\tau^2)^2.
\end{align}
where $\tau$ is sufficiently small, such that $C\tau\le \frac{1}{2}$. With the aid of Lemma  \ref{GR-KDV-lem2.3}, we have
\begin{align*}%\label{GR-KDV:eq:4.23}
||{\bm e}^n||_h^2+||\Delta_h{\bm e}^n||_{h}^2\le C(N^{-s}+\tau^2)^2,
\end{align*}
that is,
\begin{align}\label{GR-KDV:eq:4.23}
||{\bm e}^n||_h+||\Delta_h{\bm e}^n||_{h}\le C(N^{-s}+\tau^2).
\end{align}
By using Lemmas \ref{GR-KDV:lem4.3} and \ref{GR-KDV:lem4.6}, and Eq. \eqref{GR-KDV:eq:4.23}, we have
\begin{align}\label{GR-KDV:eq4.19}
||{\bm u}^{n}-{\bm U}^{n}||_{h}&+||\Delta_h({\bm u}^{n}-{\bm U}^{n})||_{h}\nonumber\\
&\le \big(||{\bm u}^{n}-({\bm u}^*)^n||_h+||\Delta_h{\bm u}^{n}-\Delta_h({\bm u}^*)^n||_h+||{\bm e}^n||_{h}+||\Delta_h{\bm e}^{n}||_{h}\big)\nonumber\\
&\le C(N^{-s}+\tau^{2}).
\end{align}
With Lemma \ref{GR-KDV-lem2.5}, we can deduce from \eqref{GR-KDV:eq4.19} that
\begin{align*}
||{\bm u}^{n}-{\bm U}^{n}||_{h,\infty}\le C(N^{-s}+\tau^{2}).
\end{align*}
This completes the proof.\qed
\end{prf}
\section{Numerical examples}\label{2SG:Sec5}
In this section, we will investigate the numerical behaviors of the LCN-MP scheme \eqref{GR-KDV-LEPSI}-\eqref{GR-KDV-LEPSI1} for the GR-KdV equation in 1D and 2D, respectively.
Also, the results are compared with some existing conservative finite difference schemes. For the LCN-MP scheme \eqref{GR-KDV-LEPSI}-\eqref{GR-KDV-LEPSI1}, we use the following iteration method to solve the linear equation:
\begin{align*}
\big({\bm I}+\mathbb{A}^2&+\frac{\tau}{2}\mathbb{B}+\frac{\tau}{2}\mathbb{L}_h\big){\bm U}^{n+\frac{1}{2},l+1}=\big({\bm I}+\mathbb{A}^2\big){\bm U}^{n}\nonumber\\
&-\frac{\tau}{2(p+2)}\Big[\text{diag}((\hat{\bm U}^{n+\frac{1}{2}})^p)\mathbb{L}_h{\bm U}^{n+\frac{1}{2},l}+\mathbb{L}_h\big[\text{diag}((\hat{\bm U}^{n+\frac{1}{2}})^p)\cdot{\bm U}^{n+\frac{1}{2},l}\big] \Big].
\end{align*}
 We take the initial iteration vector ${\bm U}^{n+\frac{1}{2},(0)}={\bm U}^{n}$ and each iteration will terminate if the infinity norm of the error between two adjacent iterative
steps is less than $10^{-14}$. Further, for a fixed iteration step $l$, the fast solver presented in Ref. \cite{JCWL17} is applied to solve the linear equations efficiently.

In order to quantify the numerical solution, we use the $e_{\infty}^{N,\tau}(t=t_n)$ to represent the $L^{\infty}$-norm of the error between the numerical solution $U_{j_1,j_2}^n$ and the exact solution
$u(x_{j_1},y_{j_2},t_n)$ at $t=t_n$. In what follows spatial mesh steps are uniformly chosen as
$h_1=h_2=h$, i.e., $N_1=N_2=N$ for simplicity.

\subsection{One dimensional R-KdV equation}
In this section, we consider the following R-KdV equation in 1D \cite{zuo09}
\begin{align}\label{RKDV:eq:5.2}
 u_t+u_{xxxxt}+u_{xxx}+u_x+uu_x=0,\ x\in\Omega,
\end{align}
with the initial condition
\begin{align*}
u(x,0)=\Big(-\frac{35}{24}+\frac{35}{312}\sqrt{313}\Big)\text{sech}^4
\Big[\frac{1}{24}\sqrt{-26+2\sqrt{313}}x\Big],\ x\in\Omega,
\end{align*}
and the periodic boundary condition.

Eq. \eqref{RKDV:eq:5.2} possesses the following exact solution \cite{zuo09}
\begin{align*}
u(x,t)=\Big(-\frac{35}{24}+\frac{35}{312}\sqrt{313}\Big)\text{sech}^4
\Big[\frac{1}{24}\sqrt{-26+2\sqrt{313}}\Big(x-\big(\frac{1}{2}+\frac{1}{26}\sqrt{313}\big)t\Big)\Big].
\end{align*}

% in this subsection, we use the $e_{2}^{h,\tau}(t=t_n)$ and $e_{\infty}^{h,\tau}(t=t_n)$ repent $L^2$- and $L^{\infty}$-norms of the error between the numerical solution $U_j^n$ and the exact solution
%$u(x_j,t_n)$ at $t=t_n$. %respectively, as
%\begin{align*}
%e_{2}^{h,\tau}(t=t_n)=h\sum_{j=0}^{N-1}|u_{j}^n-u(x_j,t_n)|^2,\ e_{h,\infty}(t_n)=\max\limits_{0\le j\le N-1}|u_{j}^n-u(x_j,t_n)|,\ 0\le n\le M.
%\end{align*}
%For the convergence rate, we use the formula
%\begin{align*}
%\text{Rate}=\frac{\ln( error_{1}/error_{2})}{\ln (\tau_{1}/\tau_{2})},
%\end{align*}
%where $\tau_{l}, error_{l}, (l=1,2)$ are step sizes and errors with the step size $\tau_{l}$, respectively.
In our computation, we take the computational domain $\Omega=[-50,50]$.
Table \ref{Tab_GR-KDV:1} shows numerical error and convergence rate of the proposed scheme with $N=1024$ and different time steps at $t=1$. As illustrated in Table \ref{Tab_GR-KDV:1}, the LCN-MP scheme has second-order
convergence rate in time. In Table \ref{Tab_GR-KDV:2}, we display the spatial numerical error and convergence rate of the proposed scheme with $\tau=10^{-5}$ and different mesh points at $t=1$, which implies that the scheme has spectral accuracy in space. { We should note that,
after $N=64$, the spatial error of the LCN-MP scheme does not decrease and is dominated by the time
discretization error. This is due to the fact that, for sufficiently smooth
problems, the Fourier pseudo-spectral method is of arbitrary order in space.} The numerical error and CPU time for different scheme with different mesh points and time steps at $t=1$ are shown in Table \ref{Tab_GR-KDV:3}. Compared with the linearized and conservative finite difference (LC-FD) scheme presented in Ref. \cite{HXH13}, our scheme provides smaller numerical error. Further, it is clear to see that, for a given $L^{\infty}$-error, the LCN-MP scheme is computationally cheaper
than the LC-FD scheme.

\begin{figure}[H]
\centering\begin{minipage}[t]{60mm}
\includegraphics[width=65mm]{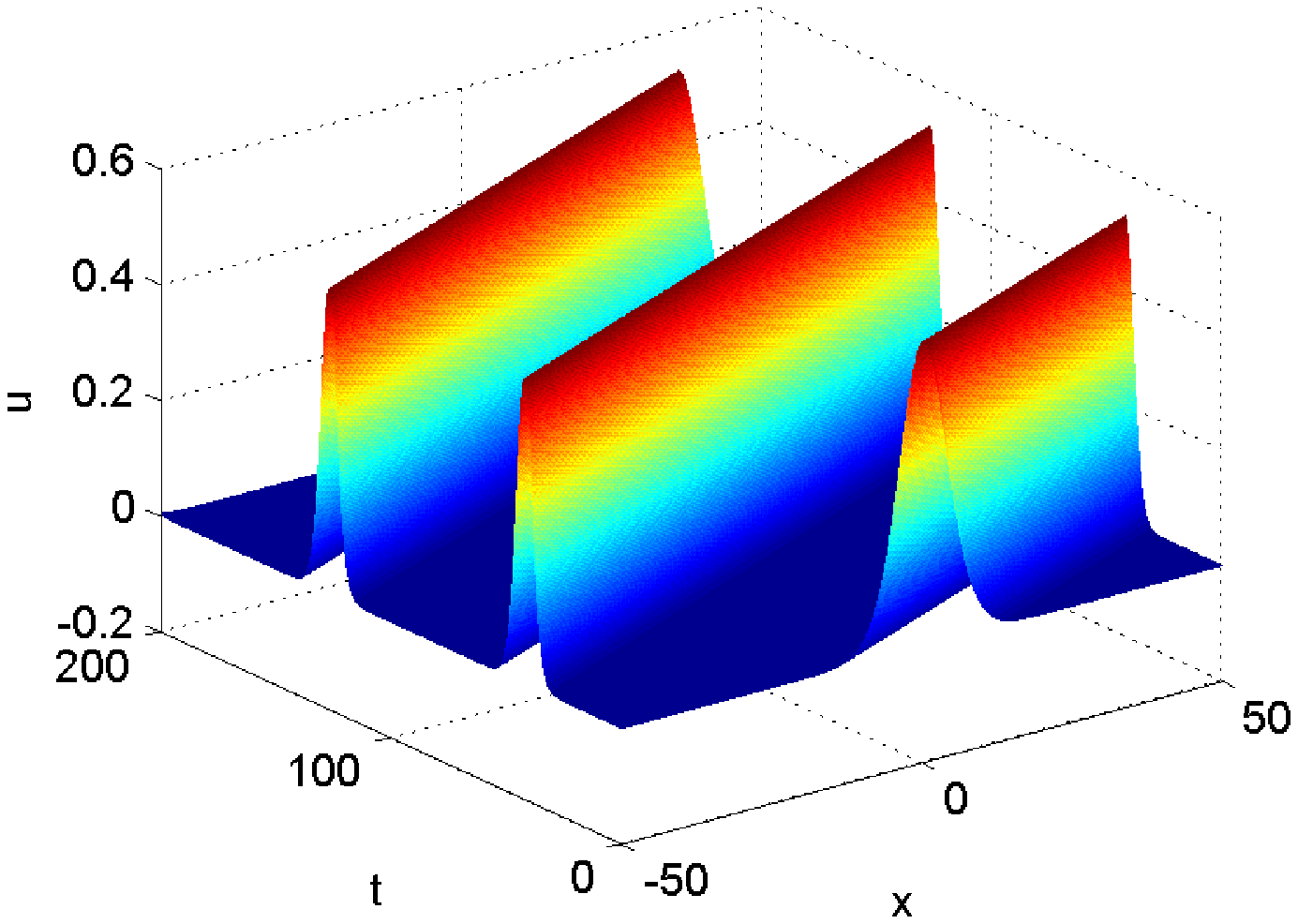}
\caption*{(a) Numerical solution}
\end{minipage}\ \
\begin{minipage}[t]{60mm}
\includegraphics[width=65mm]{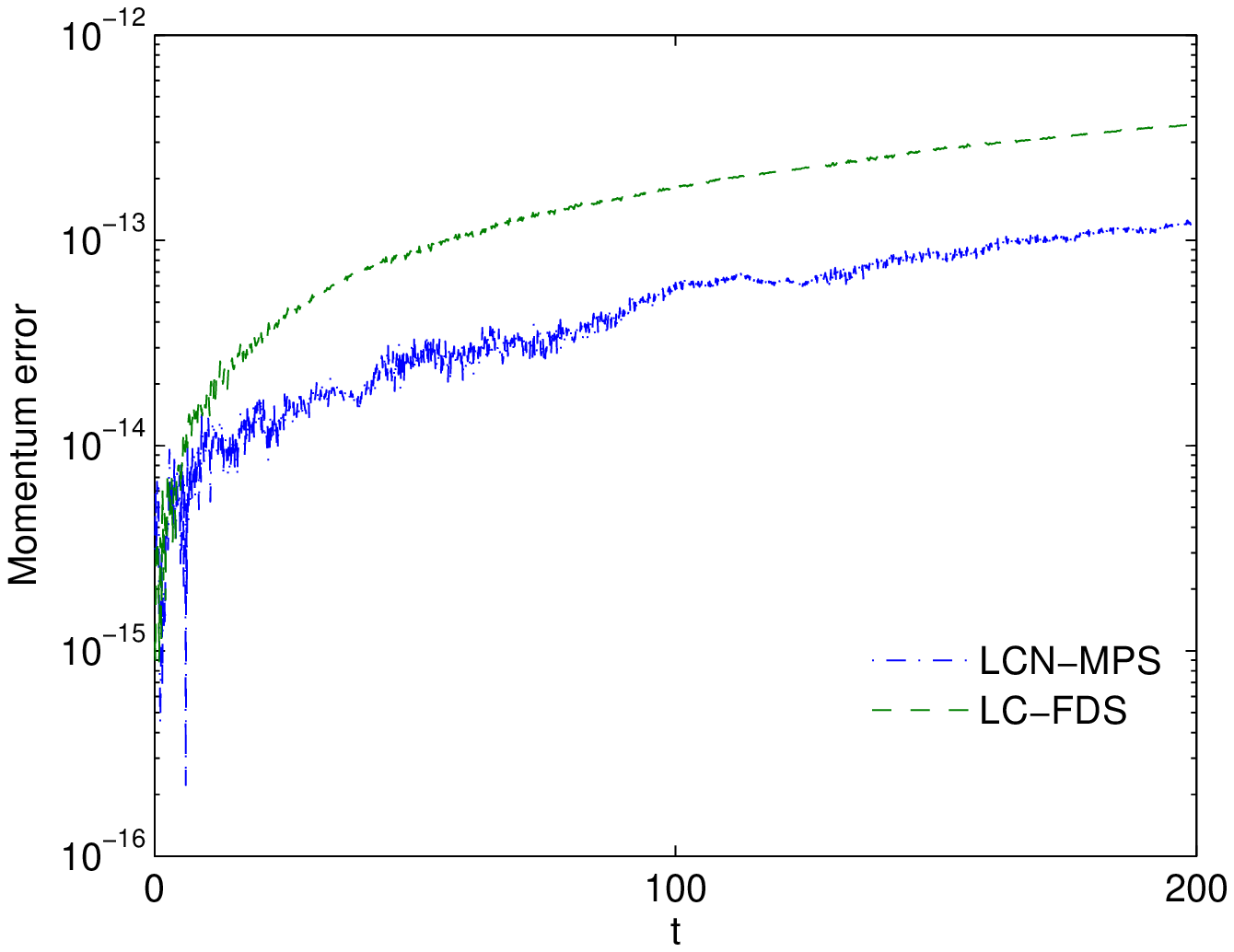}
\caption*{(b) Momentum error }
\end{minipage}
\caption{Numerical solution computed by the LCN-MP scheme (left) and momentum error (right) with $h=\tau=0.1$ over the time interval $t\in[0,200]$.}\label{Fig_GR-KDV:2}
\end{figure}

In Fig. \ref{Fig_GR-KDV:2} (a), we display the propagation of the soliton by the LCN-MP scheme over the
time interval $t\in[0,200]$, which shows that shapes of the soliton is preserved accurately in long time computation. {Here, the soliton propagates back to the computational domain because of the periodic boundary condition. Actually, for many realistic cases, perfectly matched layers \cite{Berenger94} or absorbing (or artificial) boundary conditions \cite{AABES08} have to be imposed so that the soliton can propagate throughout the computational domain. However, constructing structure-preserving schemes for the R-KdV equation under such boundary conditions is much more complied and will be our
future work.}
 The momentum error over the time interval $t\in[0,200]$ is investigated in Fig. \ref{Fig_GR-KDV:2} (b). As
illustrated in the figure, the momentum error provided by our scheme is much smaller than the one provided by the LC-FD scheme.

\begin{table}[H]
\tabcolsep=9pt
\small
\renewcommand\arraystretch{1.2}
\centering
\caption{{The temporal numerical error and convergence rate of the proposed scheme with $N=1024$ and different time steps.}}\label{Tab_GR-KDV:1}
\begin{tabular*}{\textwidth}[h]{@{\extracolsep{\fill}} c c c c c c}\hline
{$\tau$}  &{0.1} & {0.05}  & {0.025} & {0.0125} \\     %% 第1 行
\hline
 {$e_{\infty}^{N,\tau}(t=1)$}  & {2.8001e-05} & {6.9585e-06} & {1.7341e-06}&{4.3281e-07}\\[1ex]
 {Rate}  & {-} & {2.01} & {2.00}&{2.00}\\\hline
\end{tabular*}
\end{table}

\begin{table}[H]
\tabcolsep=9pt
\small
\renewcommand\arraystretch{1.2}
\centering
\caption{{{The spatial numerical error and convergence rate of the proposed scheme with $\tau=10^{-5}$ and different mesh points. }}}\label{Tab_GR-KDV:2}
\begin{tabular*}{\textwidth}[h]{@{\extracolsep{\fill}} c c c c c c}\hline
{$N$}  &{16} & {32}  & {64} & {128}&{256} \\     %% 第1行
\hline
{$e_{\infty}^{N,\tau}(t=1)$}  & {1.7538e-02} & {4.2655e-04} & {2.3645e-08}&{6.1330e-011}&{6.2242e-011}\\[1ex]
 {Rate}  & {-} & {5.4} & {14.1}&{-}&{-}\\\hline
\end{tabular*}
\end{table}

\begin{table}[H]
\tabcolsep=9pt
\small
\renewcommand\arraystretch{1.2}
\centering
\caption{{The numerical error and the CPU time of different schemes with different mesh points and time steps.}}\label{Tab_GR-KDV:3}
\begin{tabular*}{\textwidth}[h]{@{\extracolsep{\fill}} c c c c c c}\hline
{Scheme\ \ }&{$(N,\tau)$}   & {$e_{\infty}^{N,\tau}(t=1)$} &{CPU (s)} \\     %% 第1行
\hline
 {LCN-MPS}&{(1000,0.01)}  &  {2.7688e-07}&{1.9}\\[1ex]
 {}&{(2000,0.005)}  &  { 6.9176e-08}&{7.7}\\[1ex]
  {}&{(4000,0.0025)}   & {1.7289e-08}&{14.1}  \\[1ex]
 %{IFD}&{(32,0.01)}  & {4.8214e-04} & {1.5235e-04}&{1.6}\\
% {}&{(64,0.005)}  & {1.1577e-04}  & {3.6735e-05}&{6.4}\\
%  {}&{(128,0.0025)} & {2.8649e-05}   & {9.0920e-06}&{88.36}  \\
 {LC-FDS \cite{HXH13}}  & {(1000,0.01)}  &  {1.8893e-05}&{1.6}\\[1ex]
 {}&{(2000,0.005)}  &  {4.7232e-06}&{4.5}\\[1ex]
  {}&{(4000,0.0025)}   & {1.1809e-06}&{9.2}  \\\hline
\end{tabular*}
\end{table}

%\begin{figure}[H]
%\centering\begin{minipage}[t]{65mm}
%\includegraphics[width=65mm]{cpu_time}
%\caption*{(a)}
%\end{minipage}\ \
%\begin{minipage}[t]{65mm}
%\includegraphics[width=65mm]{cpu_time2}
%\caption*{(b)}
%\end{minipage}
%\caption{The CPU time versus the $L^{\infty}$-error.}\label{Fig_GR-KDV:1}
%\end{figure}

\subsection{Two dimensional GR-KdV equation}

{{\bf Example 1.}} We consider the nonhomogeneous GR-KdV equation \cite{AO15}
\begin{align}\label{RKDV-eq:5.3}
 u_t+\Delta^2u_t+\Delta u_x+(1+u^p)\mathbb{L}u=g(x,y,t),\ \ (x,y)\in\Omega,\ 0<t\leq T,
\end{align}
with the initial condition
\begin{align*}%\label{eq:1.3}
u(x,y,0)=\sin(2\pi x)\sin(2\pi y),\ (x,y)\in\Omega,
\end{align*}
and the periodic boundary conditions.

When \begin{align*}
g(x,y,t)=&\sin(2\pi x)\sin(2\pi y)\exp(-t)(-64\pi^4-1)-16\pi^3\cos(2\pi x)\sin(2\pi y)\exp(-t)\\
&+2\pi\exp(-t)\sin(2\pi(x+y))(1+\sin^p(2\pi x)\sin^p(2\pi y)\exp(-pt)),
\end{align*}
 equation \eqref{RKDV-eq:5.3} possesses the analytical solution \cite{AO15}
\begin{align*}
u(x,y,t)=\sin(2\pi x)\sin(2\pi y)\exp(-t),
\end{align*}
In our computation, we take $\Omega=[0,1]^2$ and $p=2$. Temporal and spatial numerical errors and convergence rates of the LCN-MP scheme at $t=1$ are shown in Tables  \ref{Tab_GR-KDV:4} and \ref{Tab_GR-KDV:5}, respectively. It can be observed from those tables
that the LCN-MP scheme has second-order convergence rate in time and spectral accuracy in space, respectively, which confirms the theoretical analysis. The comparisons between our scheme with the linearized and conservative finite difference (LC-FD) scheme proposed in Ref. \cite{AO15} for the numerical error and CPU time are displayed in Table \ref{Tab_GR-KDV:6}, which shows that the LCN-MP scheme has the significant advantage in the accuracy and computational efficiency over the LC-FD scheme.
\begin{table}[H]
\tabcolsep=9pt
\small
\renewcommand\arraystretch{1.2}
\centering
\caption{{The temporal numerical error and convergence rate of the proposed scheme with $N=100$ and different time steps.}}\label{Tab_GR-KDV:4}
\begin{tabular*}{\textwidth}[h]{@{\extracolsep{\fill}} c c c c c c}\hline
{$\tau$} &{0.1}  & {0.05}   & {0.025} & {0.0125} \\     %% 第1 行
\hline
 {$e_{\tau,\infty}(t=1)$}  & {4.6227e-03}&{1.1709e-03}&{2.9464e-04}&{7.3903e-05}\\[1ex]
 {Rate}  &{-} &{1.98}&{1.99}&{2.00}\\\hline
\end{tabular*}
\end{table}
\begin{table}[H]
\tabcolsep=9pt
\small
\renewcommand\arraystretch{1.2}
\centering
\caption{{{The spatial numerical error and convergence rate of the proposed scheme with $\tau=10^{-5}$ and different mesh sizes.}}}\label{Tab_GR-KDV:5}
\begin{tabular*}{\textwidth}[h]{@{\extracolsep{\fill}} c c c c c c}\hline
{$N$} &{4}  & {8}   & {16} \\     %% 第1 行
\hline
 {$e_{h,\infty}(t=1)$}  & {7.9657e-05}&{5.8725e-011}&{5.2181e-011}\\[1ex]
 {Rate}   &{-}&{20.4}&{-}\\\hline
\end{tabular*}
\end{table}
\begin{table}[H]
\tabcolsep=9pt
\small
\renewcommand\arraystretch{1.2}
\centering
\caption{{The numerical error and the CPU time of different schemes with different mesh points and time steps.}}\label{Tab_GR-KDV:6}
\begin{tabular*}{\textwidth}[h]{@{\extracolsep{\fill}} c c c c c c}\hline
{Scheme\ \ }&{$(N,\tau)$}   & {$e_{\infty}^{N,\tau}(t=1)$} &{CPU (s)} \\     %% 第1行
\hline
 {LCN-MPS}&{(8,0.01)}  &  {4.6817e-05}&{0.2}\\[1ex]
 {}&{(16,0.005)}  &  {1.1719e-05}&{0.9}\\[1ex]
  {}&{(32,0.0025)}   & {2.9615e-06}&{6.4}  \\[1ex]
 {LC-FDS \cite{HXH13}}  & {(8,0.01)}  &  {6.8370e-02}&{0.2}\\[1ex]
 {}&{(16,0.005)}  &  {1.6365e-02}&{1.1}\\[1ex]
  {}&{(32,0.0025)}   & {4.0997e-03}&{8.5}  \\\hline
\end{tabular*}
\end{table}

{\bf Example 2.}
We then consider the following GR-KdV equation in 2D \cite{AO15}:
\begin{align}\label{RKDV:eq:5.3}
 u_t+\Delta^2u_t+\Delta u_x+(1+u^p)\mathbb{L}u=0,\ \ (x,y)\in\Omega,\ 0<t\leq T,
\end{align}
with the initial condition
\begin{align*}%\label{eq:1.3}
u(x,y,0)=0.1(1+\sin(3x)\sin(5 y)),\ (x,y)\in\Omega,
\end{align*}
and the periodic boundary conditions.

\begin{figure}[H]
\centering\begin{minipage}[t]{70mm}
\includegraphics[width=70mm]{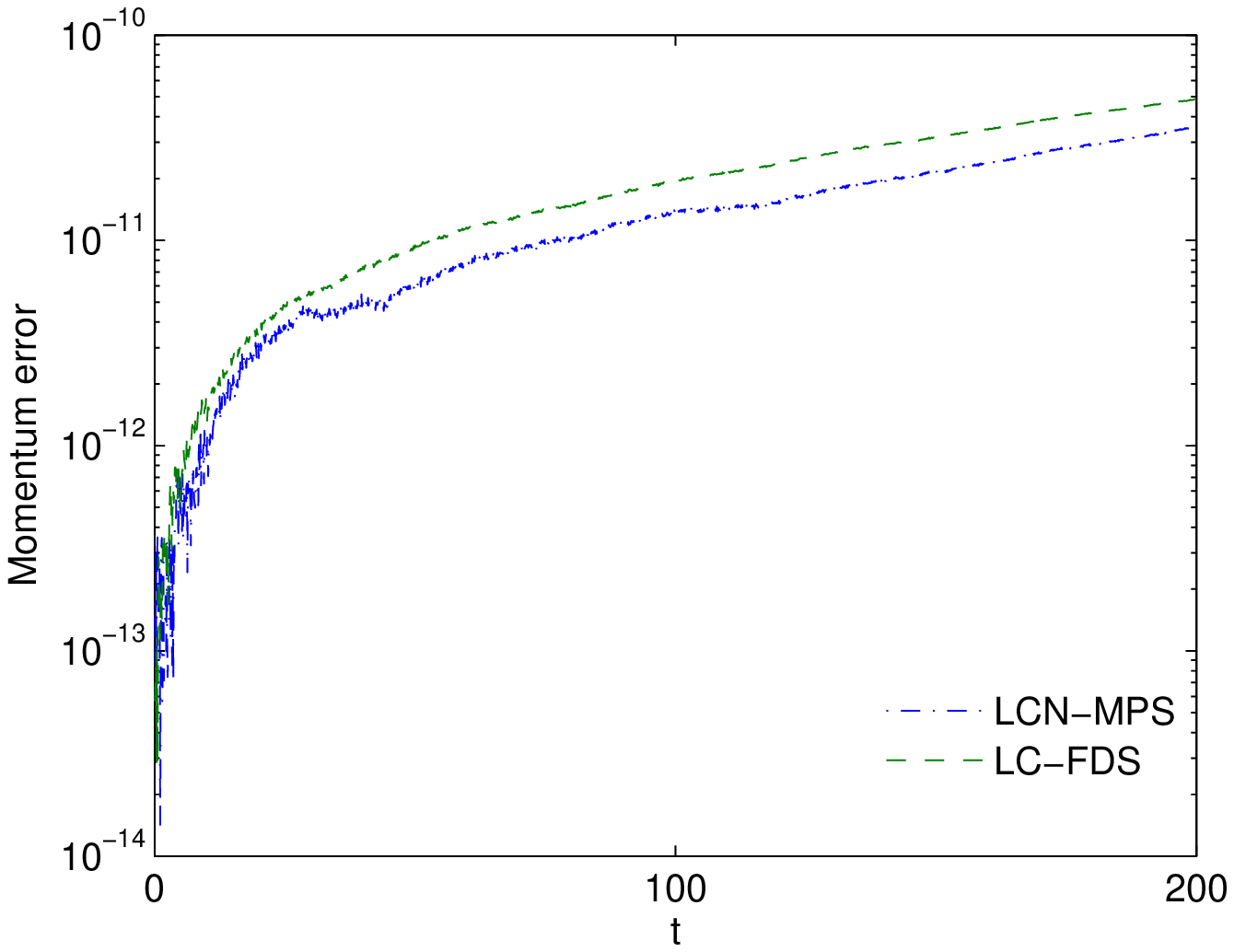}
\end{minipage}
\caption{Momentum error ($P^0=114.59$) with $h= \frac{2\pi}{50}$ and $\tau=0.1$ over the time interval $t\in[0,200]$.}\label{Fig_GKDV3}

\end{figure}

In our computation, we take $\Omega=[0,2\pi]^2$ and $p=2$. %\textcolor{blue}{The analytical solution of the GR-KdV equation \eqref{RKDV:eq:5.3} can not be obtained explicitly, thus, following Ref. \cite{ZSH14}, we define the error in the spatial direction with sufficiently $\tau$
%and the error in the temporal direction with sufficiently small $h$ at $t=t_n$, respectively, as
%\begin{align*}
%e_{h,\infty}(t=t_n)=\max\limits_{(j_1,j_2)\in J_{h}^{'}}|U_{j_1,j_2}^n(\tau,h)-U_{2j_1,2j_2}^n(\tau,h/2)|,
%\end{align*}
%and
%\begin{align*}
%e_{\tau,\infty}(t=t_n)=\max\limits_{(j_1,j_2)\in J_{h}^{'}}|U_{j_1,j_2}^n(\tau,h)-U_{j_1,j_2}^{2n}(\tau/2,h)|.
%\end{align*}
%Here, $U_{j_1,j_2}^n(\tau,h),\ (j_1,j_2)\in J_{h}^{'}$ is denoted as the numerical solution of time
%step $\tau$ and space mesh size $h$ at time $t=t_n$.} %\textcolor{blue}{ Here, we omit the comparisons between our scheme with the linearized and conservative finite difference (LC-FD) scheme proposed in Ref. \cite{AO15} for the numerical errors and CPU times. This is due to the fact that the obtained results behave similarly
%as that of Table \ref{Tab_GR-KDV:3}.}
In Fig. \ref{Fig_GKDV3}, we show the momentum error provided by the LCN-MP scheme and LC-FD scheme, respectively, over the time interval $t\in[0,200]$, which behaves similarly
as that of Fig. \ref{Fig_GR-KDV:2} (b).

\section{Concluding remarks}
In this paper, we propose a new linearized and momentum-preserving Fourier pseudo-spectral method for the GR-KdV equation. By establishing a new
semi-norm equivalence, we obtain the bound of the numerical solution in $L^{\infty}$-norm from the discrete momentum conservation law. Subsequently, based on the energy method and the bound of the numerical solution, an a priori estimate in discrete $L^{\infty}$-norm for the scheme is established without any restriction on the mesh ratio. Numerical results verify the theoretical analysis. Compared with the existing conservative schemes, our scheme is more accurate and has the significant advantage in computational efficiency and preserving the discrete momentum conservation law. Furthermore, the technique presented in this paper can also be used to establish an optimal
$L^{\infty}$-error estimate for the linearized and momentum-preserving Fourier pseudo-spectral schemes of the other Rosenau-type equation, such as the Rosenau-RLW equation \cite{PZ12}, the Rosenau-Kawahara equation \cite{BTL11,zuo09}, the Rosenau-KdV-RLW equation \cite{WD18}, etc.

%In Refs. \cite{BIT10,BI16}, Brugnano et al. have proposed a new approach, termed as Hamiltonian boundary value methods (HBVMs) to construct arbitrarily high-order energy-preserving schemes for Hamiltonian systems. Actually, when the HBVMs are extended to solve the semi-discretized system \eqref{GR-KDV_eq:2.7}, a class of arbitrarily high-order momentum-preserving schemes can be obtained. However, the numerical analysis for the resulting schemes need to be
%further discussed. Therefore, the construction and numerical analysis for high-order momentum-preserving schemes of the GR-KdV equation will be the subject of our future research.
\section*{Acknowledgments}
The authors would like to express sincere gratitude to the referees for their insightful
comments and suggestions. This work is supported by the National Natural Science Foundation of China (Grant
Nos. 11771213, 61872422), the National Key Research and Development Project of China (Grant
Nos. 2016YFC0600310, 2018YFC0603500, 2018YFC1504205), the Major Projects of Natural Sciences of
University in Jiangsu Province of China (Grant Nos. 15KJA110002, 18KJA110003), the Natural Science
Foundation of Jiangsu Province, China (Grant No. BK20171480), the Foundation of Jiangsu Key Laboratory
for Numerical Simulation of Large Scale Complex Systems (201905) and the Yunnan Provincial Department
of Education Science Research Fund Project (2019J0956).

\bibliographystyle{plain}
%\bibliography{ref}

\end{document}